\documentclass[final,authoryear,1p,times,fleqn]{elsarticle}

\usepackage{graphicx}
\usepackage{amssymb,amsthm,amsmath}
\usepackage{arydshln}

\usepackage[bookmarks=true,bookmarksnumbered=true,bookmarkstype=toc, colorlinks=true, linkcolor=black, citecolor=black, filecolor=black, pagecolor=black, urlcolor=black]{hyperref}

\newtheorem{theo}{Theorem}[section]
\newtheorem*{defi}{Definition}
\newtheorem{assu}{Assumption}[section]
\newtheorem{lemm}{Lemma}[section]

\newtheorem*{prf}{Proof}
\newtheorem*{rem}{Remarks}

\newtheorem{prop}{Proposition}[section]

\usepackage{pxfonts}
\usepackage{times}
\usepackage{bm}
\usepackage{url}
\usepackage{lineno}
\usepackage[whole]{bxcjkjatype} % pdfLatexで日本語を使う

\makeatletter
	
	\@addtoreset{equation}{section}
\makeatother

\journal{Elsevier}

\begin{document}

\begin{frontmatter}




\title{A new look at departure time choice equilibrium models with heterogeneous users\tnoteref{label1}}

\tnotetext[label1]{An earlier version of this paper has been circulated under the title ``Departure time choice equilibrium and optimal transport problems."}

\author{Takashi Akamatsu\fnref{label2}\corref{cor1}}
\ead{akamatsu@plan.civil.tohoku.ac.jp}
\author{Kentaro Wada\fnref{label3}\corref{cor1}}
\ead{wadaken@sk.tsukuba.ac.jp}
\author{Takamasa Iryo\fnref{label2}}
%\ead{takamasa.iryo.e2@tohoku.ac.jp}
\author{Shunsuke Hayashi\fnref{label4}}
%\ead{s\_hayashi@hosei.ac.jp}

\cortext[cor1]{Corresponding authors.}
\address[label2]{Graduate School of Information Sciences, Tohoku University, Aramaki Aoba 6-3-09, Aoba-ku, Sendai, Miyagi, 980-8579, Japan}
\address[label3]{Faculty of Engineering, Information and Systems, University of Tsukuba, Tennodai 1-1-1, Tsukuba, Ibaraki 305-8573, Japan}
\address[label4]{Faculty of Science and Engineering, Hosei University, Kajino-Cho 3-7-2, Koganei, Tokyo 184-8584, Japan.}

\begin{abstract}
This paper presents a systematic approach for analyzing the departure-time choice equilibrium (DTCE) problem of a single bottleneck with heterogeneous commuters. 
The approach is based on the fact that the DTCE is equivalently represented as a linear programming problem with a special structure, which can be analytically solved by exploiting the theory of optimal transport combined with a decomposition technique. 
By applying the proposed approach to several types of models with heterogeneous commuters, it is shown that (i) the essential condition for emerging equilibrium ``sorting patterns," which have been known in the literature, is that the schedule delay functions have the ``Monge property," (ii) the equilibrium problems with the Monge property can be solved analytically, and (iii) the proposed approach can be applied to a more general problem with more than two types of heterogeneities.
\end{abstract}


\end{frontmatter}

%\tableofcontents

%\newpage

%\linenumbers

%%%%%%%%%%%%%%%%%%%%%%%%%%%%%%%%%%%%%%%%%
% 1. Introduction %%%%%%%%%%%%%%%%%%%%%%%
%%%%%%%%%%%%%%%%%%%%%%%%%%%%%%%%%%%%%%%%%
\section{Introduction}\label{sec:Intro}
The modeling and analysis of rush-hour traffic congestion has a rich history, dating back to the seminal work of \cite{Vickrey1969}. 
In the basic model proposed by Vickrey, it is assumed that a fixed number of commuters with homogeneous preferences wish to arrive at a single destination (workplace) at the same preferred time, traveling through a single route that has a bottleneck of a fixed capacity. 
Each commuter chooses the departure time of his/her trip from home to minimize his/her generalized trip cost, including trip time, queuing delay at the bottleneck and schedule delay (i.e., costs of arriving early or late at their destination). 
The problem is to determine a dynamic equilibrium distribution of departure times for which no commuter can reduce his/her generalized cost by changing his/her departure time unilaterally. 
The importance of this problem in transportation planning and demand management policies has led to various extensions of the basic model to allow for (1) distributed/heterogeneous preferred departure times, (2) heterogeneity in the valuation of travel time and schedule delay, (3) elastic demands and (4) modal/route/parking choices. 
For a recent comprehensive review, refer to \cite{Li2020l}.

Despite these extensive studies, however, the extant analysis approaches have several limitations, and are not necessarily organized into a sufficiently general theory that enables us to systematically understand various extant results. 
First, most previous studies focused only on a single type of user heterogeneity either in the preferred arrival time at the destination \citep[e.g.,][]{Hendrickson1981,Smith1984a,Daganzo1985a} or in the valuation of travel time and schedule delay \citep[e.g.,][]{Cohen1987,Arnott1988,Arnott1992,Arnott1994a,Ramadurai2010,VandenBerg2011,Liu2015a,Takayama2017}\footnote{A few exceptions are \cite{Newell1987}; and \cite{Lindsey2004}; however, the former assumes a special class of schedule functions and the latter focuses only on proving the existence and uniqueness of the equilibrium.}. 
As a result, little is known about certain regularities of the equilibrium when there are multiple types of heterogeneity in the preferences of users. 
Second, most studies restricted their analysis to a special class of schedule functions (e.g., a piecewise linear function)\footnote{\cite{Smith1984a} and \cite{Daganzo1985a} are exceptions; however, they assume that all commuters have the same schedule delay function.}. 
This is theoretically problematic because it blurs the critical conditions required for the analysis results.
For instance, two models of \cite{Daganzo1985a} and \cite{Arnott1994a} have different types of users' heterogeneities (i.e., work start times and value of times (VOTs)/penalties); however, a similar sorting property emerges: 
users arrive at the bottleneck in the order of their work start times (i.e., the ``first-in-first-work (FIFW) principle") in the former model; users' arrivals are sorted in the order of their VOTs in the latter model.
This evokes a several questions, such as what conditions are essential for the emergence of ``sorting patterns" in the equilibrium and what conditions are required for obtaining an analytical solution?. 
However, answers to such basic questions are not known.
Finally, for the models with a general schedule delay function and user heterogeneity, systematic and efficient methods for obtaining accurate solutions are still lacking. 
As shown by \cite{Nie2009a}, the straight-forward formulation of the model (equilibrium conditions) results in a variational inequality (VI) problem with a non-monotone mapping, which implies that a na\"{\i}ve numerical algorithm does not guarantee convergence.

This paper proposes a systematic and unified approach to analyze a wide variety of models of departure-time choice equilibrium (DTCE) with a single bottleneck. 
The proposed approach is based on the following two facts: (1) the equilibrium can be obtained by solving a structured linear programming (LP) problem and (2) the special structure of the equivalent LP problem allows us to apply the theory of optimal transport \citep{Rachev1998,Burkard2007,Villani2008}, which provides explicit analytical solutions as well as efficient numerical algorithms. 
To attain the purpose in a clear manner, we mainly consider the DTCE models whose properties (under restricted assumptions) are familiar in the literature.

More specifically, we first review the relationships between the arrival-time-based formulation of the DTCE models and the departure-time-based one. 
It has been shown that the departure-time-based formulation has a big advantage, that is, the (part of) equilibrium conditions reduce to an equivalent LP problem; the LP problem can be effectively utilized to study pricing and system optimal problems and can be solved easily by an established convergent algorithm (or an LP solver).
We then reveal that the equivalent LP problems for various types of DTCE problems have structural commonalities, to which the theory of optimal transport can be applied either directly or indirectly. 
Through several examples, we demonstrate how the optimal transport theory can be applied to the DTCE problems.
The above-mentioned sorting patterns for different models can be explained by a property of the schedule delay function in a unified manner; even analytical solutions are obtained for certain cases.
We also show that the proposed approach enables us to systematically analyze the basic properties of the equilibrium, such as existence, uniqueness, and regularities of the flow patterns (e.g., temporal ``sorting" patterns).

As a first example of the approach, we demonstrate that the DTCE problem with heterogeneous preferred departure times is analytically solvable and that the sorting property of the equilibrium flow pattern (i.e., the FIFW principle) can be understood from the optimal transport theory as a direct consequence of ``submodularity" or the ``Monge property" of the schedule delay cost function.  
As a second example, we consider the DTCE problem with heterogeneous schedule cost functions (i.e., users are differentiated according to their VOT). 
For this type of DTCE problem, the straightforward application of the optimal transport theory is not possible because the schedule cost functions do not satisfy the Monge property. 
Nevertheless, by developing an approach in which the optimal transport theory is applied to the subproblems generated from a hierarchical decomposition of the original LP problem, we can show that the DTCE problem is analytically solvable under mild assumptions regarding the schedule cost function. 
Finally, this approach is further extended to generalized DTCE models with two types of cost heterogeneity.

In the remainder of this article, after briefly describing the problem setting (Section \ref{sec:setting}), Section \ref{sec:arrival} reviews the arrival-time-based formulation of the DTCE models. 
With some new variables, Section \ref{sec:DTCandLP} converts the arrival-time-based formulation to the departure-time-based one.
We then introduce an LP approach to the DTCE problem with heterogeneous commuters and show its some advantages and applications. 
Section \ref{sec:OTP} briefly reviews the theory of optimal transport and presents an illustrative example of how it can be applied to analyze the DTCE problem with heterogeneous preferred departure times. 
Section \ref{sec:VOT} analyzes the DTCE problems with heterogeneous schedule delay cost functions, in which we provide a new analytical approach of combining a hierarchical decomposition with the optimal transport theory. 
Section \ref{sec:3D_Monge} further extends the approach to generalized DTCE models with two types of cost heterogeneity. 
Concluding remarks are presented in Section \ref{sec:Conclusion}.

%%%%%%%%%%%%%%%%%%%%%%%%%%%%%%%%%%%%%%%%%
% 2.  Setting            %%%%%%%%%%%%%%%%
%%%%%%%%%%%%%%%%%%%%%%%%%%%%%%%%%%%%%%%%%
\section{Morning commute equilibrium problems}\label{sec:setting}
%\subsection{Setting of the model}
We consider a road network with a single O-D pair connected by one route. 
We assume that one user travels in every vehicle. 
Users are treated as a continuum, and the total mass is denoted by $Q$, which is a given constant.
The route has a single bottleneck with a capacity (maximum service rate) of $\mu$.
The queuing congestion at the bottleneck is described by a point queue model, in which the queue is assumed to form vertically at the entrance of the bottleneck.
Users are classified into a finite number, $K$, of homogeneous groups. 
The index set of groups is denoted by $\mathcal{K} = \{1, 2, \dots, K\}$.
Let the mass of users in group $k \in \mathcal{K}$ be $Q_{k}$; then, $\sum_{k\in \mathcal{K}} Q_{k} = Q$.

Users choose their departure time from the origin (or the bottleneck) so as to minimize their trip cost.
The trip cost is assumed to be additively separable into free-flow travel, queuing delay, and schedule delay costs.
The free-flow travel time from the origin to the bottleneck is assumed, without loss of generality, to be zero unless otherwise noted.
We also assume that the departure time from the bottleneck is the arrival time at the destination. 
The schedule delay, $\epsilon$, is defined as the difference between the actual and preferred departure times from the bottleneck; the schedule delay cost function $c_{k}(\epsilon)$ measured in the (queuing) time unit is assumed to be a continuous and piecewise differentiable function of the schedule delay and is specified in the later sections.

Under these assumptions, the {\it morning commute equilibrium} is defined as the state in which no user could reduce his/her trip cost by changing his/her departure/arrival time unilaterally. 
There are two types of formulations for the equilibrium: {\it arrival-time-based} and {\it departure-time-based}. 
In the former formulation, all variables are expressed as functions of the arrival time at the bottleneck; in the latter formulation, they are functions of the departure time from the bottleneck. 
While the two formulations are essentially equivalent, their mathematical structures are quite different for some instances. 
In the next section, we will show the former formulation; we will then convert it to the latter one, which is the basis of our new approach, in Section \ref{sec:DTCandLP}.

%%%%%%%%%%%%%%%%%%%%%%%%%%%%%%%%%%%%%%%%%
% 3.  Arrival-time-base  %%%%%%%%%%%%%%%%
%%%%%%%%%%%%%%%%%%%%%%%%%%%%%%%%%%%%%%%%%
\section{Arrival-time-based formulation}\label{sec:arrival}
\subsection{Formulation}
Let $t\in \mathcal{T}$ be the {\it arrival time at the bottleneck} and $\mathcal{T} \subset \mathbb{R}$ be a sufficiently long morning rush hour during which all users complete their trips. 
The queuing delay of users with arrival time $t$ at the bottleneck is denoted by $d(t)$. 
The schedule delay $\epsilon$ of users in group $k$ whose the preferred departure time $\sigma_{k}$ is expressed as $t + d(t) - \sigma_{k}$. 
Because $\sigma_{k}$ is constant for users, their schedule delay cost function is a function of the departure time, $t + d(t)$; that is, $c_{k}(\epsilon) = c_{k}(t + d(t))$.

At equilibrium, the following four conditions should hold.
The first condition is the uses' optimal choice condition:
\begin{align}
	& \begin{cases}
		v_{k} = d(t) + c_{k}(t + d(t)) \mspace{15mu} \text{if} \mspace{15mu} \lambda_{k}(t) > 0\\
		v_{k} \leq d(t) + c_{k}(t + d(t)) \mspace{15mu} \text{if} \mspace{15mu} \lambda_{k}(t) = 0
	\end{cases} & \forall k\in \mathcal{K}, \ \forall t\in \mathcal{T}, \label{eq:Arr_UserChoice}
\end{align}
where $v_{k}$ represents the minimum (equilibrium) trip cost for users in group $k$ and $\lambda_{k}(t)$ denotes group $k$'s arrival flow rate at the bottleneck at time $t$.

The second condition is the queuing delay condition. 
In the point queue model (for the details, see \ref{sec:PQ}), the queuing delay $d(t)$ for a user arriving at the bottleneck at time $t$ can be represented by
\begin{subequations}
\begin{align}
	&\dot{d}(t) = \begin{cases}
		\lambda(t)/\mu - 1 \mspace{85mu} \text{if} \mspace{15mu} d(t) > 0\\
		\max. [0, \lambda(t)/\mu - 1] \mspace{15mu} \text{if} \mspace{15mu} d(t) = 0
	\end{cases} & \forall t\in \mathcal{T}, 
	\label{eq:OriginalPQ}
\end{align}
where the ``dot" denotes the derivative operation with respect to arrival time $t$, and $\lambda(t)\equiv \sum_{k\in \mathcal{K}}\lambda_{k}(t)$. 
This implies that the queuing delay should satisfy
\begin{align}
	&\begin{cases}
		\dot{d}(t) = \lambda(t)/\mu - 1 \mspace{15mu} \text{if} \mspace{15mu} d(t) > 0\\
		\dot{d}(t) \geq \lambda(t)/\mu - 1 \mspace{15mu} \text{if} \mspace{15mu} d(t) = 0
	\end{cases} & \forall t\in \mathcal{T}.
	\label{eq:RelaxPQ}
\end{align}
\end{subequations}
In this study, we employ the complementarity condition \eqref{eq:RelaxPQ}, instead of \eqref{eq:OriginalPQ}, as the queuing delay model.
This is because it is analytically tractable and its essential features are consistent with the original point queue model (e.g., the first-in-first-out (FIFO) property holds)\footnote{See \cite{Ban2012}, \cite{Han2013} and \cite{Jin2014f} for more detailed discussions.}.

The third condition is the flow conservation condition. 
For each user group $k$, the integral of the arrival flow rate $\lambda_{k}(t)$ at the bottleneck must be equal to the total mass $Q_{k}$; that is,
\begin{align}
	& \int_{\mathcal{T}}\lambda_{k}(t)\mathrm{d}t = Q_{k} & \forall k\in \mathcal{K}.
	\label{eq:Arr_conservation}
\end{align}

The last condition determines the departure flow rate, $\kappa(t)$, from the bottleneck.
This condition seems trivial, but it is necessary to complete the equilibrium flow pattern. 
Following the FIFO principle, we have the following flow propagation condition:
\begin{subequations} \label{eq:Construct_Dep}
\begin{align} 
	& D(t + d(t)) = A(t) & \forall t \in \mathcal{T}\\
	& \kappa(t + d(t)) = \lambda(t)/(1 + \dot{d}(t)) & \forall t \in \mathcal{T},
\end{align} 
\end{subequations} 
where $A(t)$ [$D(t)$] denotes the cumulative arrival [departure] flow at the bottleneck by time $t$, $\lambda(t) = \mathrm{d}A(t)/\mathrm{d}t$, and $\kappa(t) = \mathrm{d}D(t)/\mathrm{d}t$.  
Note that the resulting $\kappa(t)$ is always nonnegative and satisfies the capacity constraint (i.e., $\kappa(t) \leq \mu$).

\begin{figure}[t]
\begin{center}
\includegraphics[width = 100mm,clip]{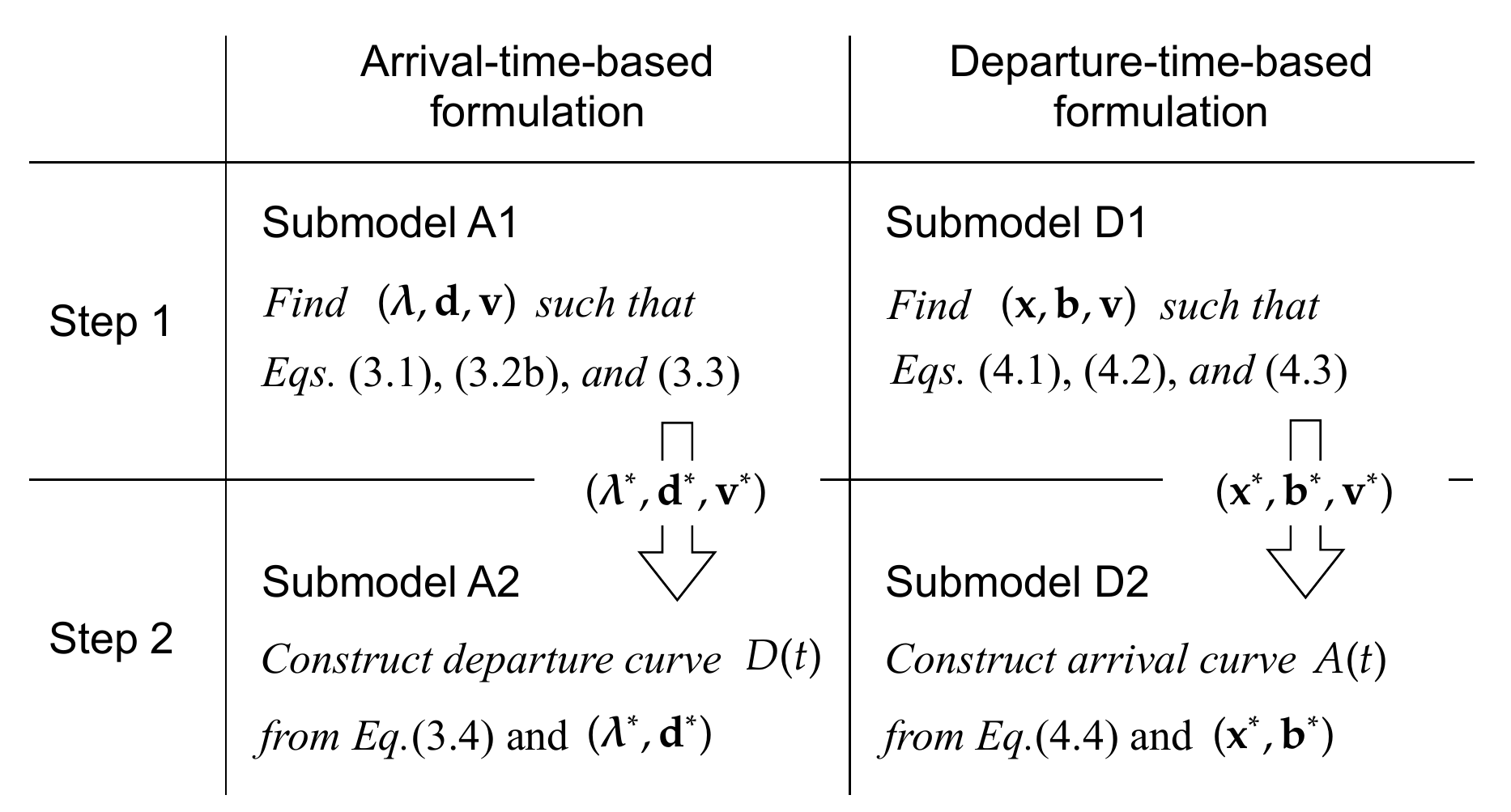}
\caption{Correspondence between two equilibrium formulations}\label{fig:Framework}
\vspace{-4mm}
\end{center}
\end{figure}

Because there are only three types of the unknown variables $(\bm{\lambda}, {\bf d}, {\bf v})$, the first three conditions are sufficient to determine them.
Therefore, we can solve the equilibrium problem in two sequential steps (see Figure \ref{fig:Framework}). 
The first step (Submodel A1) solves the equilibrium conditions \eqref{eq:Arr_UserChoice}, \eqref{eq:RelaxPQ}, and \eqref{eq:Arr_conservation} simultaneously. 
The second step (Submodel A2) constructs the equilibrium departure pattern by combining the condition \eqref{eq:Construct_Dep} with the solution from the first step. 
Alternatively, we can obtain $\kappa(t)$, or equivalently $D(t)$, from the point queue model for a given $A(t)$ in a more constructive manner.

\subsection{Equivalent complementarity and variational inequality problems of Submodel A1}
Submodel A1, \eqref{eq:Arr_UserChoice}, \eqref{eq:RelaxPQ} and \eqref{eq:Arr_conservation}, can be expressed as the following continuous-time (or infinite-dimensional) nonlinear complementarity problem (NCP): {\it Find $(\bm{\lambda}, {\bf d}, {\bf v})$ such that}
\begin{subequations}\label{eq:NCP}
	\begin{align}
		& 0 \leq \lambda_{k}(t) \ \bot \ d(t) + c_{k}(t + d(t)) - v_{k}\geq 0 & \forall k\in \mathcal{K}, \ \forall t\in \mathcal{T},\\
		& 0\leq d(t) \ \bot \ \dot{d}(t) - [\lambda(t)/\mu - 1] \geq 0 & \forall t\in \mathcal{T},\\
		& 0 \leq v_{k} \ \bot \ \textstyle{\int_{\mathcal{T}}\lambda_{k}(t)\mathrm{d}t - Q_{k}} \geq 0 & \forall k\in \mathcal{K}.
	\end{align}
\end{subequations}
It can be seen that, if the schedule delay cost function $c_{k}(t + d(t))$ is linear in the schedule delay (and thus linear in the queuing delay $d(t)$), the NCP reduces to the linear complementarity problem (LCP).  
For example, with a piecewise linear schedule delay cost function, the discrete-time problem is formulated as a finite-dimensional LCP \citep{Ramadurai2010}; its continuous-time counterpart was studied as a linear complimentarily system \citep{Pang2011}.
However, there is no study on the NCP \eqref{eq:NCP}, and its properties are largely unknown.

Submodel A1 can also be studied as a continuous-time variational inequality (VI) problem as follows. 
%Following \cite{Nie2009a}, a VI formulation of the morning commute problem is given as follows. 
{\it Find $\bm{\lambda}^{*} \in \Omega_{\lambda}$ such that}
\begin{align}
	& \sum_{k\in \mathcal{K}}\int_{\mathcal{T}}C_{k}(\bm{\lambda})\left(\lambda_{k}(t) - \lambda_{k}^{*}(t)\right)\mathrm{d}t \geq 0 & \forall \bm{\lambda}\in \Omega_{\lambda}, 
\end{align}
where $C_{k}(\bm{\lambda})$ is the trip cost of users in group $k$, which is an implicit function of the arrival flow rate $\bm{\lambda}$ and is evaluated by the schedule delay cost function and queuing delay condition, and $\Omega_{\lambda} = \{\lambda_{k}(t) \ | \ \textstyle{\int_{\mathcal{T}}\lambda_{k}(t)\mathrm{d}t = Q_{k}}, \ \lambda_{k}(t)\geq 0, \ \forall k\in \mathcal{K}, \ \forall t \in \mathcal{T}\}$. 
However, because the mapping $C_{k}(\bm{\lambda})$ of this VI problem may not be monotone even for the simplest case with $K = 1$ (i.e., homogeneous user case) and a piecewise linear schedule delay cost function \citep{Nie2009a}, it is generally difficult to establish a convergent solution method for it.

The above observations indicate that the arrival-time-based formulation of the morning commute problem may not be a tractable way to analyze the equilibrium for a general setting. 
In contrast, the departure-time-based formulation shown in the next section possesses a much simpler mathematical expression of the (part of) equilibrium conditions, that is, it reduces to an equivalent LP problem.

%%%%%%%%%%%%%%%%%%%%%%%%%%%%%%%%%%%%%%%%%
% 4.  Departure-time-base  %%%%%%%%%%%%%%
%%%%%%%%%%%%%%%%%%%%%%%%%%%%%%%%%%%%%%%%%
\section{Departure-time-base formulation with linear programming}\label{sec:DTCandLP}
\subsection{Formulation}
Let $s\in \mathcal{S}$ be the {\it departure time from the bottleneck} and $\mathcal{S} \subset \mathbb{R}$ be a sufficiently long morning rush hour during which all users complete their trips. 
The queuing delay of users with departure time $s$ from the bottleneck is denoted by $b(s)$.
As in the previous section, because the schedule delay $\epsilon$ is $s-\sigma_{k}$, the schedule delay cost function is expressed as a function of the departure time, $s$; that is, $c_{k}(\epsilon) = c_{k}(s)$.
However, unlike the arrival-time-based formulation, $c_{k}(s)$ is evaluated independent of the state variables and thus can be treated as a given {\it parameter} (or a given function). 
We will convert the arrival-time-based formulation \eqref{eq:Arr_UserChoice}, \eqref{eq:RelaxPQ}, \eqref{eq:Arr_conservation}, and \eqref{eq:Construct_Dep} to the departure-time based formulation.

\begin{figure}[t]
\begin{center}
\includegraphics[width = 68mm,clip]{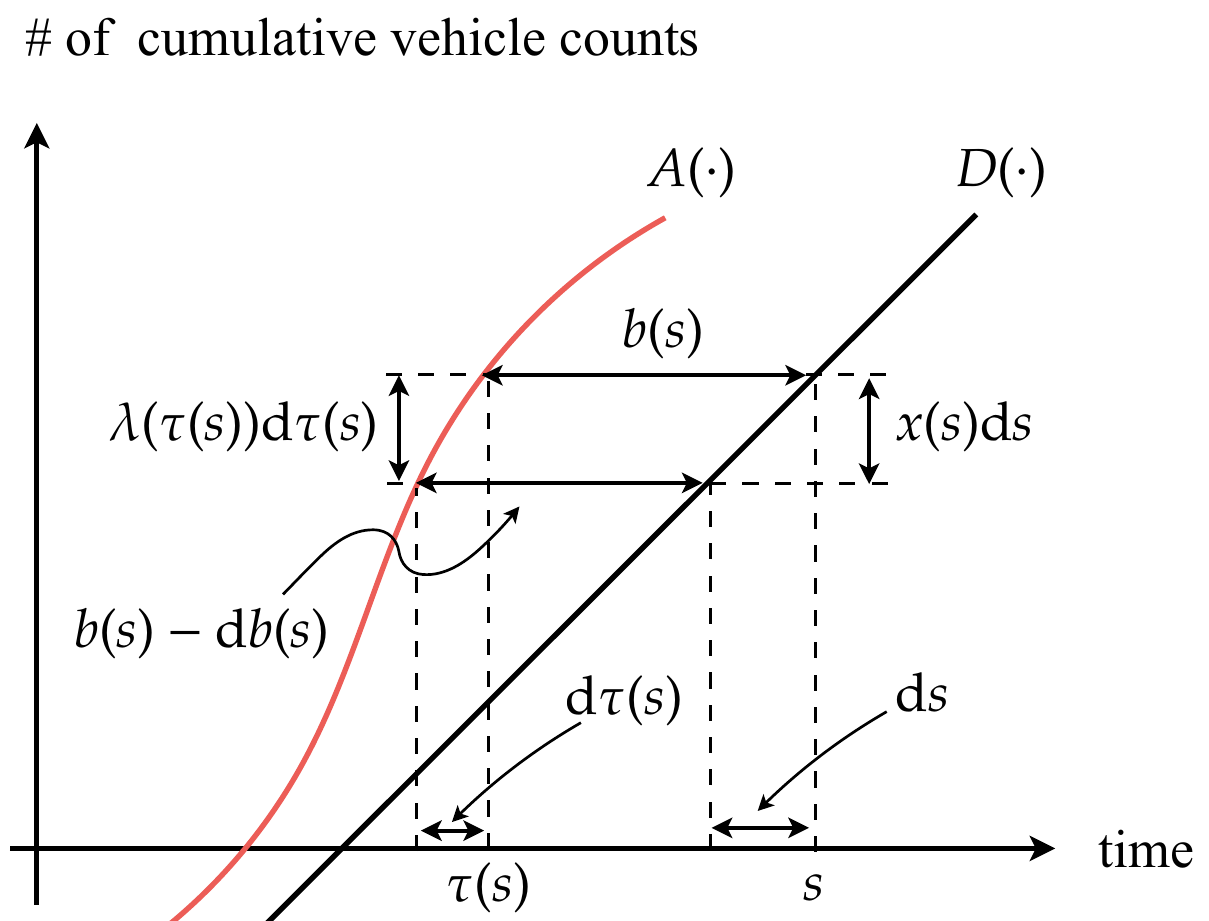}
\caption{Relationships between arrival and departure times/flows}\label{fig:Relation}
\vspace{-4mm}
\end{center}
\end{figure}

As a preliminary step in the conversion, we summarize the relationships between the variables in the arrival-time-based and departure-time-based formulations. 
First, the relationship between the arrival and departure times is given by $t = \tau(s) \equiv s - b(s)$. 
Then, the FIFO condition can be written as, $D(s) = A(\tau(s))$. 
By taking a derivative of both sides of this condition with respect to $s$, we have
\begin{align*}
	& x(s) = \lambda(\tau(s))\cdot \Delta \tau(s), 
\end{align*}
where $x(s) \equiv \Delta D(s)$ denotes the departure flow rate from the bottleneck at time $s$ and $\Delta$ denotes the derivative operation with respect to {\it bottleneck-departure-time} $s$; that is,  
\begin{align*}
	& \Delta \tau(s) \equiv \mathrm{d}\tau(s)/\mathrm{d}s = 1 - \Delta b(s). 
\end{align*}
These relationships can be graphically understood from Figure \ref{fig:Relation}. 
Regarding the cost variables, for users departing from the bottleneck at time $s$, the queuing delay $b(s)$ is obviously given by $b(s) = d(\tau(s))$, which implies that 
\begin{align*}
	& \Delta b(s) = \dot{d}(\tau(s)) \cdot \Delta \tau (s). 
\end{align*}

Now, let us derive the departure-time-based formulation by substituting the above relationships into the arrival-time-based formulation.
First, the users' optimal choice condition \eqref{eq:Arr_UserChoice} reduces to
\begin{subequations}\label{eq:UserChoice}
\begin{align}
	& \begin{cases}
		v_{k} = b(s) + c_{k}(s) \mspace{15mu} \text{if} \mspace{15mu} x_{k}(s)/\Delta \tau(s) > 0\\
		v_{k} \leq b(s) + c_{k}(s) \mspace{15mu} \text{if} \mspace{15mu} x_{k}(s)/\Delta \tau(s) = 0
	\end{cases} & \forall k\in \mathcal{K}, \ \forall s\in \mathcal{S}. \label{eq:UserChoiceA}
\end{align}
At equilibrium, the following condition (i.e., a relaxed version of the condition \eqref{eq:UserChoiceA}) must also be met because the departure flow rate $x_{k}(s)$ is nonnegative.
\begin{align}
	& \begin{cases}
		v_{k} = b(s) + c_{k}(s) \mspace{15mu} \text{if} \mspace{15mu} x_{k}(s) > 0\\
		v_{k} \leq b(s) + c_{k}(s) \mspace{15mu} \text{if} \mspace{15mu} x_{k}(s) = 0
	\end{cases} & \forall k\in \mathcal{K}, \ \forall s\in \mathcal{S}. \label{eq:UserChoiceB}
\end{align}
\end{subequations}
With the assumption $\Delta \tau(s) > 0$, both conditions are consistent. 
Otherwise, the condition \eqref{eq:UserChoiceA} [\eqref{eq:UserChoiceB}] can produce negative departure [arrival] flow rates. 
This means that $\Delta \tau(s) > 0$ must be met if the equilibrium exists.
We will later discuss the condition under which the assumption holds, but we accept it for the moment and use the condition \eqref{eq:UserChoiceB} as the users' optimal choice condition.

Second, the queuing delay condition \eqref{eq:RelaxPQ} can be represented as 
\begin{subequations}
\begin{align} 
	&\begin{cases}
		\Delta b(s) = x(s)/\mu - \Delta \tau (s) \mspace{15mu} \text{if} \mspace{15mu} b(s) > 0\\
		\Delta b(s)\geq x(s)/\mu - \Delta \tau (s) \mspace{15mu} \text{if} \mspace{15mu} b(s) = 0
	\end{cases} & \forall s\in \mathcal{S}.
\end{align}
Because $\Delta \tau(s) + \Delta b(s) = 1$ by definition, the condition reduces to 
\begin{align}
	& \begin{cases}
		\sum_{k\in \mathcal{K}}x_{k}(s) = \mu \mspace{15mu} \text{if} \mspace{15mu} b(s) > 0\\
		\sum_{k\in \mathcal{K}}x_{k}(s) \leq \mu \mspace{15mu} \text{if} \mspace{15mu} b(s) = 0
	\end{cases} & \forall s \in \mathcal{S}. \label{eq:PQ-capacity}
\end{align}
\end{subequations}

Third, it is easily seen that the flow conservation condition \eqref{eq:Arr_conservation} reduces to
\begin{align}
	& \int_{\mathcal{S}}x_{k}(s)\mathrm{d}s = Q_{k} & \forall k\in \mathcal{K}. \label{eq:conservation}
\end{align}

Finally, we must determine the arrival flow rate $\lambda(\tau(s))$ at the bottleneck. 
As shown above, the FIFO discipline implies that 
\begin{subequations} \label{eq:Construct_Arr}
\begin{align} 
	& A(s-b(s)) = D(s) & \forall s \in \mathcal{S}\\
	& \lambda(s - b(s)) = x(s)/(1 - \Delta b(s)) & \forall s \in \mathcal{S}.
\end{align} 
\end{subequations}
As noted, to ensure the nonnegativity and finiteness of the arrival flow rate, the above assumption, $\Delta \tau(s) = 1 - \Delta b(s) > 0$, must be met.

\begin{figure}[t]
\begin{center}
\includegraphics[width = 65mm,clip]{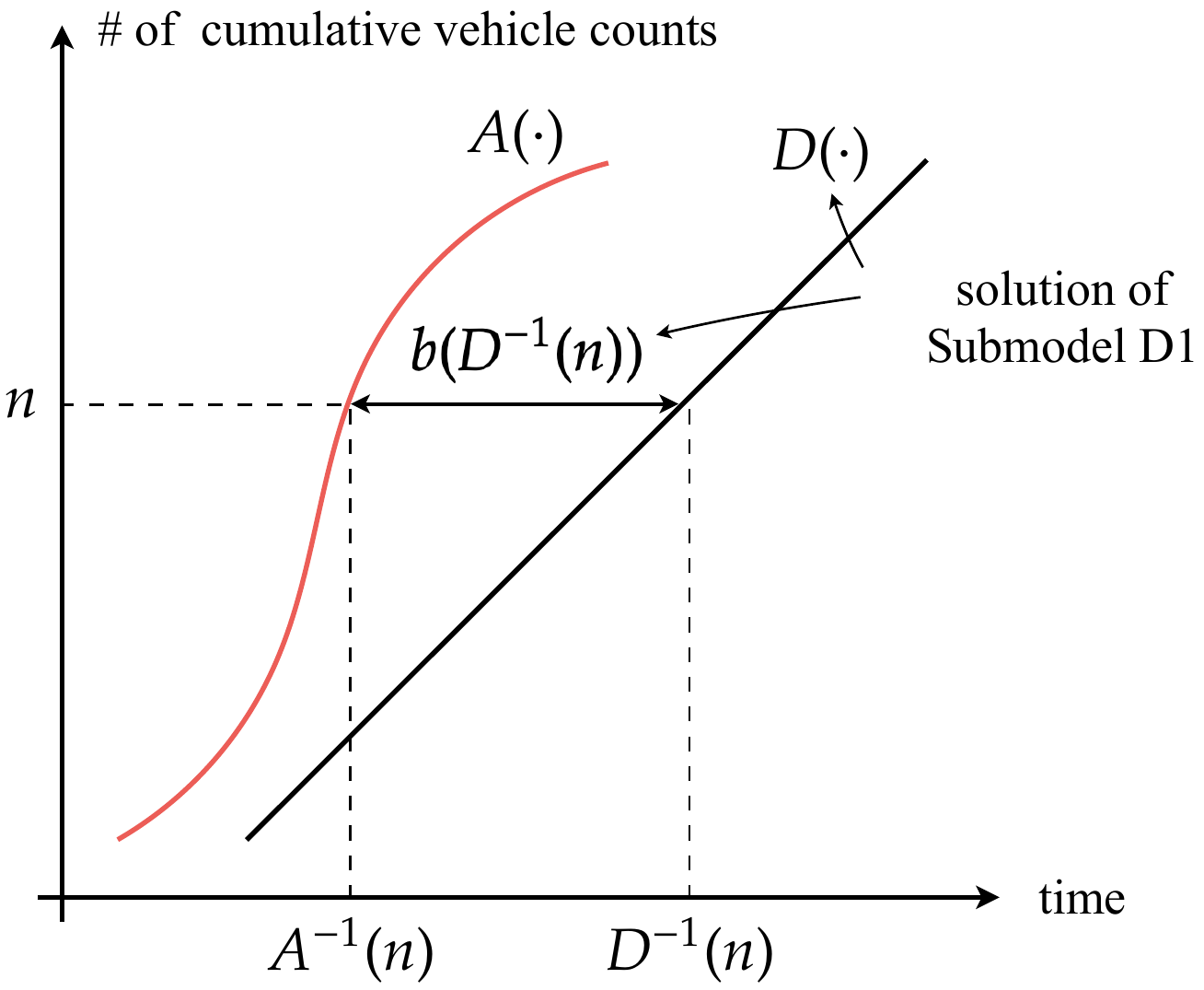}
\caption{Schematic illustration of the construction of the cumulative arrival curve from the solution of Subproblem D1}\label{fig:ArivalCurve_Construction2}
\vspace{-4mm}
\end{center}
\end{figure}

Similar to the arrival-time-based formulation, the problem can be solved in two sequential steps (see again Figure \ref{fig:Framework}). 
The first step (Submodel D1) solves the first three conditions \eqref{eq:UserChoiceB}, \eqref{eq:PQ-capacity}, and \eqref{eq:conservation}, simultaneously to determine the three types of unknown variables $({\bf x},{\bf b}, {\bf v})$. 
The second step (Submodel D2) completes the equilibrium arrival pattern by combining the condition \eqref{eq:Construct_Arr} with the solution in the first step. 
Note that a more constructive way to obtain the arrival flow pattern is to use the inverse functions of cumulative curves $A(\cdot)$ and $D(\cdot)$.
As shown in Figure \ref{fig:ArivalCurve_Construction2}, the FIFO condition of the point queue model is represented as follows:
\begin{align}
	& A^{-1}(n) = D^{-1}(n) - b(D^{-1}(n)) & \forall n, \label{eq:FIFO}
\end{align}
where $A^{-1}(n)$ and $D^{-1}(n)$ denote the $n^{\text{th}}$ user's arrival and departure times at/from the bottleneck, respectively.
Because we know both terms on the right-hand side of Eq.\eqref{eq:FIFO} (i.e., the solution in the first step), we can construct the equilibrium cumulative arrival curve $A(\cdot)$, which is perfectly consistent with the queuing model.

So far, we have shown the correspondence of each of the four equilibrium conditions with the two formulations by assuming $\Delta \tau(s) > 0$.
Now, let us show the condition under which the assumption holds.

\begin{lemm}\label{lemm:existence}
	Suppose $\Delta c_{k}(s) > -1$. Then, $\Delta \tau(s) > 0$ holds true. 
\end{lemm}

\begin{prf} {\rm See \ref{sec:proof_deltau}.}
\end{prf}

We should note that the condition, $\Delta c_{k}(s) > -1$, is consistent with the existence conditions of the equilibrium shown in \cite{Smith1984a} and \cite{Lindsey2004}.
Then, we obtain the following proposition on the equivalence of the two formulations.
\begin{prop} \label{prop:equivalence}
Suppose there exists a morning commute equilibrium, or equivalently $\Delta c_{k}(s) > -1$. 
Then, the departure-time-based equilibrium conditions \eqref{eq:UserChoiceB}, \eqref{eq:PQ-capacity}, \eqref{eq:conservation}, \eqref{eq:Construct_Arr} and the arrival-time-based conditions \eqref{eq:Arr_UserChoice}, \eqref{eq:RelaxPQ}, \eqref{eq:Arr_conservation}, \eqref{eq:Construct_Dep} are equivalent in the sense that the solution of one formulation satisfies all conditions of the other formulation.   
\end{prop}

\subsection{Equivalent linear programming of Submodel D1}
In the following sections, we propose a systematic approach to analyze a wide variety of models of a morning commute equilibrium. 
The basis of this approach is the fact that Submodel D1, \eqref{eq:UserChoiceB}, \eqref{eq:PQ-capacity}, and \eqref{eq:conservation}, reduces to an equivalent LP problem. 
The equivalent LP problem was first shown by \cite{Iryo2005} and \cite{Iryo2007} in a discrete-time setting.
We briefly describe the LP problem below, but for a {\it continuous-time} setting.

Consider first the following infinite-dimensional LP problem:
\begin{align}
\text{[2D-LP(${\bf x}$)]} \mspace{30mu} & \min_{{\bf x}\geq {\bf 0}}. \ Z({\bf x})\equiv  \sum_{k\in \mathcal{K}}\int_{\mathcal{S}}c_{k}(s)x_{k}(s)\mathrm{d}s \label{eq:2DP-obj}\\
& \text{subject to} \mspace{15mu} \sum_{k\in \mathcal{K}}x_{k}(s) \leq \mu & \forall s\in \mathcal{S} \label{eq:2DP-con1}\\
&\mspace{80mu} \int_{\mathcal{S}}x_{k}(s)\mathrm{d}s = Q_{k} & \forall k\in \mathcal{K} \label{eq:2DP-con2},
\intertext{and the associated dual problem:}
\text{[2D-LP(${\bf u}, {\bf v}$)]} \mspace{30mu} & \max_{{\bf u}\geq {\bf 0}, {\bf v}}. \ \check{Z}({\bf u}, {\bf v}) \equiv -\mu \int_{\mathcal{S}}u(s)\mathrm{d}s + \sum_{k\in \mathcal{K}}Q_{k}v_{k}\\
& \text{subject to} \mspace{15mu}  c_{k}(s) + u(s) - v_{k} \geq 0 & \forall k\in \mathcal{K}, \ \forall s\in \mathcal{S}, \label{eq:2DPD-con1}
\end{align}
where $u(s)$ and $v_{k}$ are Lagrange multipliers for \eqref{eq:2DP-con1} and \eqref{eq:2DP-con2}, respectively. 
As noted, $c_{k}(s)$ is a parameter (or a coefficient) and is given exogenously.   
As shown in \ref{sec:Duality}, the strong duality for [2D-LP(${\bf x}$)] holds, implying the following complementarity slackness (or optimality) conditions: 
\begin{align}
	& \begin{cases}
		x_{k}(s) \{c_{k}(s) + u(s) - v_{k}\} = 0\\
		\mspace{55mu} c_{k}(s) + u(s) - v_{k}\geq 0, \ x_{k}(s) \geq 0
	\end{cases} & \forall k\in \mathcal{K}, \ \forall s\in \mathcal{S} \label{eq:KT-1}\\
	& \begin{cases}
		u(s) \{\mu - \textstyle{\sum_{k\in \mathcal{K}}}x_{k}(s)\} = 0\\
		\mspace{55mu} \mu - \textstyle{\sum_{k\in \mathcal{K}}}x_{k}(s) \geq 0, \ u(s)\geq 0
	\end{cases} & \forall s\in \mathcal{S}. \label{eq:KT-2}
\end{align}

Comparing the optimality conditions above and the equilibrium conditions \eqref{eq:UserChoiceB}, \eqref{eq:PQ-capacity}, and \eqref{eq:conservation}, we can easily see that these are equivalent if we can interpret the Lagrange multiplier ${\bf u}$ of the optimality conditions as the queuing delay {\bf b}. 
For this interpretation to be valid, the arrival flow rate at the bottleneck should be physically feasible (i.e., non-negative and finite).
From the same logic as in the proof of {\bf Lemma \ref{lemm:existence}}, it is easily seen that the interpretation is valid if the schedule cost function satisfies $\Delta c_{k}(s) > -1$\footnote{$\Delta c_{k}(s) > -1$ asserts $\Delta \tau(s) = 1 - \Delta u(s) > 0$, where $\tau(s) = s - u(s)$, and thus $\lambda(\tau(s)) = x(s)/\Delta \tau(s) \geq 0$ is obtained.}. 
The above discussion is summarized in the following proposition.

\begin{prop} \label{prop:equivalence1}\citep{Iryo2005,Iryo2007,Akamatsu2015a}
Suppose there exists a morning commute equilibrium, or equivalently $\Delta c_{k}(s) > -1$. 
Then, the optimal solution $({\bf x}, {\bf u}, {\bf v})$ of [2D-LP(${\bf x}$)] and [2D-LP(${\bf u}, {\bf v}$)] is consistent with the equilibrium conditions.
\end{prop}

It is now clear that the departure-time-based formulation has an advantage over the arrival-based formulation. 
Submodel D1 with {\it a general schedule delay function} reduces to [2D-LP(${\bf x}$)], which can be solved and analyzed much easier than the equivalent NCP or VI of Submodel A1. 
In Sections \ref{sec:OTP} and \ref{sec:VOT}, we will study specific morning commute equilibrium problems by appropriately setting $c_{k}(s)$ of the equivalent LP problem [2D-LP(${\bf x}$)] for each problem. 
However, it should be emphasized here that [2D-LP(${\bf x}$)] has a powerful ability to deal with a wide variety of models (other than the models in Sections \ref{sec:OTP} and \ref{sec:VOT}) without imposing restrictive assumptions on schedule delay functions (or user heterogeneities).

\subsection{Equilibrium versus system optimal}
The equivalent LP problem is more powerful than it appears; it can also be applied to some control problems. 
One important class of such examples is the equilibrium problem under tolls \citep[e.g.,][]{Arnott1990,Laih1994,Lindsey2012,Chen2015a}.
Let $p(s)$ be the exogenously given toll for departure time $s$ from the bottleneck. 
Then, the users' optimal choice condition at equilibrium is given by
\begin{align}
	& \begin{cases}
		x_{k}(s) \{c_{k}(s) + \alpha_{k}^{-1}p(s) + b(s) - v_{k}\} = 0\\
		\mspace{55mu} c_{k}(s) + \alpha_{k}^{-1}p(s) + b(s) - v_{k}\geq 0, \ x_{k}(s) \geq 0
	\end{cases} & \forall k\in \mathcal{K}, \ \forall s\in \mathcal{S}, \label{eq:UserChoiceB-toll}
\end{align}
where $\alpha_{k}$ denotes the VOT for users in group $k$. 
It is easy to see that the (departure-time-based) equilibrium formulation in which Eq.\eqref{eq:UserChoiceB} is replaced with Eq.\eqref{eq:UserChoiceB-toll} has the following equivalent LP formulation:  
\begin{align}
	& \min_{{\bf x}\geq {\bf 0}}. \ Z_{\text{Toll}}({\bf x})\equiv  \sum_{k\in \mathcal{K}}\int_{\mathcal{S}}\left\{c_{k}(s) + \alpha_{k}^{-1}p(s)\right\}x_{k}(s)\mathrm{d}s \label{eq:EP_toll}\\
	& \text{subject to} \mspace{15mu} \text{Eq.\eqref{eq:2DP-con1} and Eq.\eqref{eq:2DP-con2}}. \notag
\end{align}
Note that, as in the case of [2D-LP(${\bf x}$)], $\Delta c_{k}(s) + \alpha_{k}^{-1}\Delta p(s) > -1$ must be met if we interpret the Lagrangian multiplier $u(s)$ for the capacity constraint \eqref{eq:2DP-con1} as the queuing delay $b(s)$. 
While the above-mentioned studies do not use the LP formulation, there is a possibility of obtaining more insight into the theoretical properties of the problem (e.g., some regularity of traffic flow patterns) by exploiting the LP approach shown in the following sections.

Instead of the {\it time}-based equilibrium formulation above, we can describe the equilibrium under tolls by using costs measured in {\it monetary} unit. 
Specifically, the users' optimal choice condition in the {\it monetary}-based formulation can be written as 
\begin{align}
	& \begin{cases}
		x_{k}(s) \{\alpha_{k}c_{k}(s) + p(s) + \alpha_{k}b(s) - \alpha_{k}v_{k}\} = 0\\
		\mspace{55mu} \alpha_{k}c_{k}(s) + p(s) + \alpha_{k}b(s) - \alpha_{k}v_{k}\geq 0, \ x_{k}(s) \geq 0
	\end{cases} & \forall k\in \mathcal{K}, \ \forall s\in \mathcal{S}. \label{eq:UserChoiceB-toll-money}
\end{align}
This is equivalent to the condition \eqref{eq:UserChoiceB-toll} because users' optimal choice does not depend on the scaling of costs for each user group\footnote{
It is thus clear that the equilibrium formulation (or [2D-LP(${\bf x}$)]) in previous sections includes the cases in which users are differentiated by the combination of VOTs and schedule delay functions. 
}. 
Note, however, that the {\it monetary}-based formulation has no equivalent optimization problem, as in Eq.\eqref{eq:EP_toll} (see Section \ref{sec:comparison} for a similar discussion).

The LP representation of the morning commute equilibrium problem is also closely related to a dynamic system optimal (DSO) problem with heterogeneous users. 
Specifically, the DSO problem ({\it without queuing}) is formulated as:  
\begin{align}
	& \min_{{\bf x}\geq {\bf 0}}. \ Z_{\text{SO}}({\bf x})\equiv  \sum_{k\in \mathcal{K}}\int_{\mathcal{S}}\left\{\alpha_{k}c_{k}(s)\right\}x_{k}(s)\mathrm{d}s \label{eq:DSO}\\
	& \text{subject to} \mspace{15mu} \text{Eq.\eqref{eq:2DP-con1} and Eq.\eqref{eq:2DP-con2}}. \notag
\end{align}
The only difference from [2D-LP(${\bf x}$)] is that the objective function is measured in {\it monetary} unit, which results in that the equilibrium and optimal flow patterns being different in general. 
However, [2D-LP(${\bf x}$)] can be interpreted as a particular case of the DSO problem: the solutions of the two problems are identical if the VOTs are all the same among user groups (i.e., $\alpha_{k} = \alpha$).

The optimality conditions for the DSO problem are given by Eq.\eqref{eq:2DP-con2}, Eq.\eqref{eq:KT-2}, and 
\begin{align}
	& \begin{cases}
	x_{k}(s) \{\alpha_{k}c_{k}(s) + \tilde{u}(s) - \tilde{v}_{k}\} = 0\\
	\mspace{55mu} \alpha_{k}c_{k}(s) + \tilde{u}(s) - \tilde{v}_{k}\geq 0, \ x_{k}(s) \geq 0
\end{cases} & \forall k\in \mathcal{K}, \ \forall s\in \mathcal{S},  \label{eq:DSO-optimality}
\end{align} 
where $(\tilde{{\bf u}}, \tilde{{\bf v}})$ are the Lagrange multipliers for the constraints \eqref{eq:2DP-con1} and \eqref{eq:2DP-con2}, respectively. 
From a comparison to Eq.\eqref{eq:UserChoiceB-toll-money} (or Eq.\eqref{eq:UserChoiceB-toll}), it is readily seen that the condition \eqref{eq:DSO-optimality} is equivalent to the particular case of the users' optimal choice condition: $\tilde{u}(s) = p(s)$ and $b(s) = 0$, $\forall s$ (and $\tilde{v}_{k} = \alpha_{k}v_{k}$, $\forall k$). 
This means that $\tilde{u}(s)$ can be regarded as the dynamic optimal toll pattern $p(s)$ that eliminates the queuing delay at equilibrium \citep{Doan2011}.
Another interpretation is the market clearing price pattern under a time-dependent tradable bottleneck permit scheme, which is designed to resolve the problem of congestion during the morning rush hour at a single bottleneck \citep{Akamatsu2006,Akamatsu2007,Wada2010,Akamatsu2017}.  
Specifically, the optimality condition of the DSO problem, which corresponds to Eq.\eqref{eq:KT-2}, is interpreted as the demand-supply market equilibrium for each time $s$, where the demand of time period $s$ permit is equal to the departure flow $x(s) = \sum_{k\in \mathcal{K}}x_{k}(s)$, and the maximum supply of the permit is given by the bottleneck capacity $\mu$.

%%%%%%%%%%%%%%%%%%%%%%%%%%%%%%%%%%%%%%%%%
% 5.  Optimal Transport Problems  %%%%%%%
%%%%%%%%%%%%%%%%%%%%%%%%%%%%%%%%%%%%%%%%%
\section{Monge-Kantorovich problem}\label{sec:OTP}
In the previous section, we expressed the (part of) morning commute equilibrium problem as the linear program [2D-LP(${\bf x}$)]. 
If one only wishes to obtain a numerical solution of the equilibrium problem, it is enough to solve [2D-LP(${\bf x}$)] using an appropriate LP solver.
However, in order to obtain a deeper insight into the theoretical properties (such as uniqueness, some regularity of flow patterns) of [2D-LP(${\bf x}$)] (or the equilibrium), the theory of optimal transport \citep[see,][]{Kantorovich1942,Kantorovich1948,Rachev1998,Burkard2007,Villani2008} is useful.
In this section, we briefly review the role of Monge properties in optimization in Section \ref{sec:Monge} and present an illustrative example of its application to the departure time choice equilibrium (DTCE) in Section \ref{sec:FIFW}.

\subsection{Monge property and analytical solution}\label{sec:Monge}
The optimal transport problem in a two-dimensional discrete space setting is the following finite-dimensional LP problem, which is well known as ``{\it Hitchcock's transportation problem}" in the operations research and transportation fields:
\begin{align}
\text{[2D-OTP]} \mspace{30mu} & \min_{{\bf x}\geq {\bf 0}}. \ Z_{2D}({\bf x})\equiv  \sum_{i\in \mathcal{I}}\sum_{k\in \mathcal{K}}c_{i,k}x_{i,k} \label{eq:2D-TP_obj}\\
& \text{subject to} \mspace{15mu} \sum_{k\in \mathcal{K}}x_{i, k} = S_{i} & \forall i\in \mathcal{I} = \{1, 2, \dots, I\} \\
&\mspace{82mu} \sum_{i\in \mathcal{I}}x_{i, k} = Q_{k} & \forall k\in \mathcal{K} = \{1, 2, \dots, K\}
\end{align}
where vectors ${\bf S}$ and ${\bf Q}$ are given constants satisfying $\sum_{i\in \mathcal{I}}S_{i} = \sum_{k\in \mathcal{K}}Q_{k}$.

Before reviewing some useful theorems on the transportation problem, we introduce the following concepts:

\begin{defi}\rm An $I \times K$ real matrix $C = [c_{i, k}]$ is called a {\it Monge matrix} if $C$ satisfies the following property ({\it Monge property})
\begin{align}
	& c_{i, k} + c_{i + 1, k + 1} \leq c_{i, k+1} + c_{i+1, k} & \text{\rm{for all }} \ 1\leq i < I, \ 1 \leq k < K \label{eq:MongeProperty}
\end{align}
In addition, if the inequality in \eqref{eq:MongeProperty} strictly holds, $C$ is termed a {\it strict Monge matrix}.
\end{defi}

\begin{defi}\rm Function $c: \mathbb{R}^{2} \rightarrow \mathbb{R}$ is {\it submodular} if, and only if,
\begin{align}
	& c(x, y) + c(x', y') \leq c(x, y') + c(x', y) & \text{\rm{for all}} \ x \leq x', y \leq y' \label{eq:Submodular} 
\end{align}
In addition, if the inequality in \eqref{eq:Submodular} strictly holds for $x < x', y < y'$, $c$ is called a {\it strict submodular function}.
\end{defi}

\noindent This implies that an $I \times K$ matrix $C$ whose elements are given by $c_{i,k}:= c(i\Delta x, k\Delta y) \ (1\leq i \leq I, \ 1\leq k \leq K)$ is a (strict) Monge matrix if the function $c: \mathbb{R}^{2} \rightarrow \mathbb{R}$ is (strict) submodular\footnote{Throughout the paper, we express the variable (or parameter) as a function if it depends on a {\it continuous} index (i.e., time), while we use the subscript for describing its dependence on a {\it discrete} index (i.e., user group).}.
Thus, we will also term condition \eqref{eq:Submodular} the (continuous) Monge property.
If the inequalities \eqref{eq:MongeProperty} and \eqref{eq:Submodular} hold in the opposite direction, then the matrix $C$ and function $c$ are said to be an {\it inverse} Monge matrix and {\it supermodular}, respectively.

As is well known, a feasible solution to the transportation problem [2D-OTP] can always be determined by a greedy algorithm termed the {\it northwest corner rule} \citep{Hoffman1963}. 
\begin{description}
	\item[Northwest corner rule] \hspace{-0.4em} 
\begin{enumerate}
	\setcounter{enumi}{-1}
	\item Initialize the indices with $i := 1$ and $k := 1$.
	\item Set $x_{i, k} := \min \{S_{i}, Q_{k}\}$. 
	\item Reduce both the supply $S_{i}$ and demand $Q_{k}$ by $x_{i, k}$: $S_{i} := S_{i} - x_{i, k}$ and $Q_{k} := Q_{k} - x_{i, k}$. If some of $S_{i}$ and $Q_{k}$ become zero, then these indices are increased by one.
	\item If there still exists an unsatisfied constraint, go back to Step 2.  
\end{enumerate}
\end{description}
The Monge property further provides the following useful result:
\begin{theo}\label{theo:northwest}{\rm \citep{Hoffman1963,Hoffman1985}:} If the cost matrix $C$ of [2D-OTP] has the Monge property, then the northwest corner rule yields an optimal solution for arbitrary supply ${\bf S}$ and demand ${\bf Q}$ vectors. 
\end{theo}

\begin{rem}{\rm It is worth noting that {\it the cost matrix $C$ is not used at all in the northwest corner rule}. This implies that even if we only know that $C$ is a Monge matrix, {\it we can determine an optimal solution of [2D-OTP] without knowing the explicit values of the cost coefficients}.
If the cost coefficients fulfill the {\it inverse Monge property}, an optimal solution can be found by the {\it northeast corner rule}\footnote{In the northeast corner rule, the algorithm begins with the indices with $i := 1$ and $k := K$ (or $i := I$ and $k := 1$), then the index is decreased by one if $Q_{k}$ (or $S_{i}$) becomes zero in Step 2 in the northwest corner rule.}.
}	
\end{rem}

\begin{theo}\label{theo:northwest_uniqueness}
{\rm \citep{Dubuc1999}:} If the matrix $C = [c_{i, k}]$ in Eq.\eqref{eq:2D-TP_obj} is a strict Monge matrix, then the optimal solution of [2D-OTP] is unique (i.e., the solution provided by the northwest corner rule is the only solution of [2D-OTP]). 
\end{theo}

A continuous analog of Theorem \ref{theo:northwest} for the continuous transportation problem (i.e., all the indices of the parameters and variables in the problem are continuous) is as follows.
Let $x\in X = \mathbb{R}$ and $y\in Y = \mathbb{R}$ be random variables and let $F_{1}, F_{2}$, and $F$ denote the distribution functions of $x, y,$ and $(x, y)$, respectively.
Given a continuous cost function $c: \mathbb{R}^{2}\rightarrow \mathbb{R}$, the continuous optimal transport problem can be formulated as: 
\begin{align}
	\text{[2D-COTP]} \mspace{30mu} & \min_{F\in \mathcal{F}(F_{1}, F_{2})}. \ Z_{2D}(F)\equiv \int_{X\times Y} c(x, y)\mathrm{d}F(x, y) \label{eq:2D-CTP_obj}\\
& \text{where} \mspace{15mu} \mathcal{F}(F_{1}, F_{2}) \equiv \left\{\textstyle F(x, \infty) = F_{1}(x), \  F(\infty, y)  = F_{2}(y), \forall x, y\in \mathbb{R}\right\} \label{eq:2D-CTP_con}
\end{align}
Note that $F_{1}(\infty) = F_{2}(\infty)$. 
We can easily see the correspondence between the problems [2D-COTP] and [2D-OTP] by changing the discrete indices to the continuous ones and by rewriting the demand and supply constraints in [2D-OTP] in a {\it cumulative} form (see also Subsection \ref{sec:FIFW} for an example).

\begin{theo}\label{theo:Frechet-Hoeffding}{\rm \citep{Cambanis1976,Dubuc1999}} Let $F_{1}$ and $F_{2}$ be distribution functions on $\mathbb{R}$. Furthermore, suppose that $\sup_{F\in \mathcal{F}(F_{1}, F_{2})} Z_{2D}(F) < \infty$ and the cost function $c: \mathbb{R}^{2}\rightarrow \mathbb{R}$ is a submodular.
Then an optimal solution $F^{*}\in \mathcal{F}(F_{1}, F_{2})$ is given by the so-called Fr\'echet--Hoeffding distribution
\begin{align}
	& F^{*}(x, y) = \min \{F_{1}(x), F_{2}(y)\} & \forall (x, y)\in \mathbb{R}^{2}.
\end{align}
Furthermore, if the function $c$ is a strict submodular, the solution is unique (i.e., $F^{*}(x, y)$ is the only optimal solution of [2D-COTP]).
\end{theo}
We note that the northwest corner rule for the discrete transportation problem can be viewed as explicit rules for calculating the analytical solution (i.e., Fr\'echet--Hoeffding distribution) in the discrete space setting \citep{Burkard2007}.

\subsection{Illustrative example: Model with heterogeneous preferred departure time} \label{sec:FIFW} 
To illustrate the usefulness of the theory of optimal transport in analyzing the DTCE problem, we here consider a model with heterogeneous preferred departure times.
A well-known sorting property of the equilibrium flow pattern, the so-called the {\it FIFW principle} \citep{Daganzo1985a}, can be understood from the theory of optimal transport as a direct consequence of the Monge property of the schedule cost function.

\begin{assu} \label{ass:schedule_FIFW} The schedule delay cost function is assumed to be identical for all users and is given by
\begin{align}
	c_{k}(s) \equiv f(\epsilon(\sigma_{k}, s)), \mspace{15mu} \epsilon(\sigma_{k}, s)\equiv s - \sigma_{k}, \label{eq:schedule_FIFW}
\end{align}
where $\sigma_{k}$ is the preferred departure time of users in group $k$ and a function $f: \mathbb{R}\rightarrow \mathbb{R}$ is continuous, strictly convex, and has a minimum at $\epsilon = 0$. 
\end{assu}

\begin{lemm} \label{lemm:FIFW_Monge} Suppose that user groups are arranged (indexed) in the increasing order of their desired arrival times: $\sigma_{1} < \cdots < \sigma_{K}$. 
The function $c_{k}(s)$ defined by \eqref{eq:schedule_FIFW} satisfies the strict Monge property.
\end{lemm}
\begin{prf} \rm See \ref{sec:proof_l3.1}.     
\end{prf}

In the above setting, in which the departure time is continuous, but the preferred departure time is discrete \citep{Newell1987,Lindsey2004}, all groups experience a queuing delay to equilibrate the trip costs of users within each group.
This implies that several disjoint departure time windows can exist in equilibrium.
However, we consider a typical morning rush-hour problem in which a single joint departure time window (i.e., a single rush period) $\hat{\mathcal{S}}\equiv [s_{0}, s_{0} + T]\subset \mathcal{S}$ of length $T = Q/\mu$ occurs in equilibrium \citep[e.g.,][]{Smith1984a,Daganzo1985a}.

For a given equilibrium rush period $\hat{\mathcal{S}}$, the DTCE problem  (i.e., [2D-LP(${\bf x}$)]) can be reduced to an instance of the optimal transport problem, where $c_{k}(s)$ satisfies the {\it strict Monge property}.
\begin{align}
& \min_{{\bf x}\geq {\bf 0}}. \ Z({\bf x})\equiv  \sum_{k\in \mathcal{K}}\int_{s\in \hat{\mathcal{S}}}c_{k}(s)x_{k}(s)\mathrm{d}s \label{eq:OTP-FIFW-obj}\\
& \text{subject to} \mspace{15mu} \sum_{k\in \mathcal{K}}x_{k}(s) = \mu & \forall s\in \hat{\mathcal{S}}  \\
& \mspace{80mu}\int_{s\in \hat{\mathcal{S}}}x_{k}(s)\mathrm{d}s = Q_{k} & \forall k\in \mathcal{K}
\end{align}
\noindent Its cumulative form is expressed as follows: 
\begin{align}
& \min_{F\in \mathcal{F}}\int_{[1,K]\times \hat{\mathcal{S}}}c(\kappa, s)\mathrm{d}F(\kappa, s) \label{eq:Cum_OTP-FIFW-obj}\\
& \text{where} \mspace{15mu} \mathcal{F}\equiv \left\{F(\kappa, T) = {\textstyle \sum_{k' \leq \kappa}}Q_{k'}, \ F(K, s)  = \mu \cdot (s - s_{0}), \ \forall \kappa \in [1, K], \ \forall s\in \hat{\mathcal{S}}\right\}.
\end{align}
Note that the parameter $c_{k}(s)$ and variable $x_{k}(s)$ that depend on both discrete and continuous indices $(k, s)$ are related to $c(\kappa, s)$ and $F(\kappa, s)$ that depend on the continuous indices $(\kappa, s)$ as follows: $F(\kappa, s) \equiv \sum_{k' \leq \kappa}\int_{s_{0}}^{s}x_{k'}(s)\mathrm{d}s$, which is a step function of $\kappa$ for a given $s$; without loss of generality, $c(\kappa, s)$ is assumed to be a continuous and strict submodular function that satisfies $c(k, s) = c_{k}(s), \ \forall k = 1, \dots, K, \ \forall s\in \hat{\mathcal{S}}$\footnote{Because $\mathrm{d}F(\kappa, s) = 0$ for interval $\kappa\in (k, k+1), \forall k = 1, \dots, K-1$, by definition, the function $c(\kappa, s)$ can be set arbitrary except for the values at the points $\kappa = k, \ \forall k = 1, \dots , K$.}.

\begin{figure}[t]
\begin{center}
\includegraphics[width = 82mm,clip]{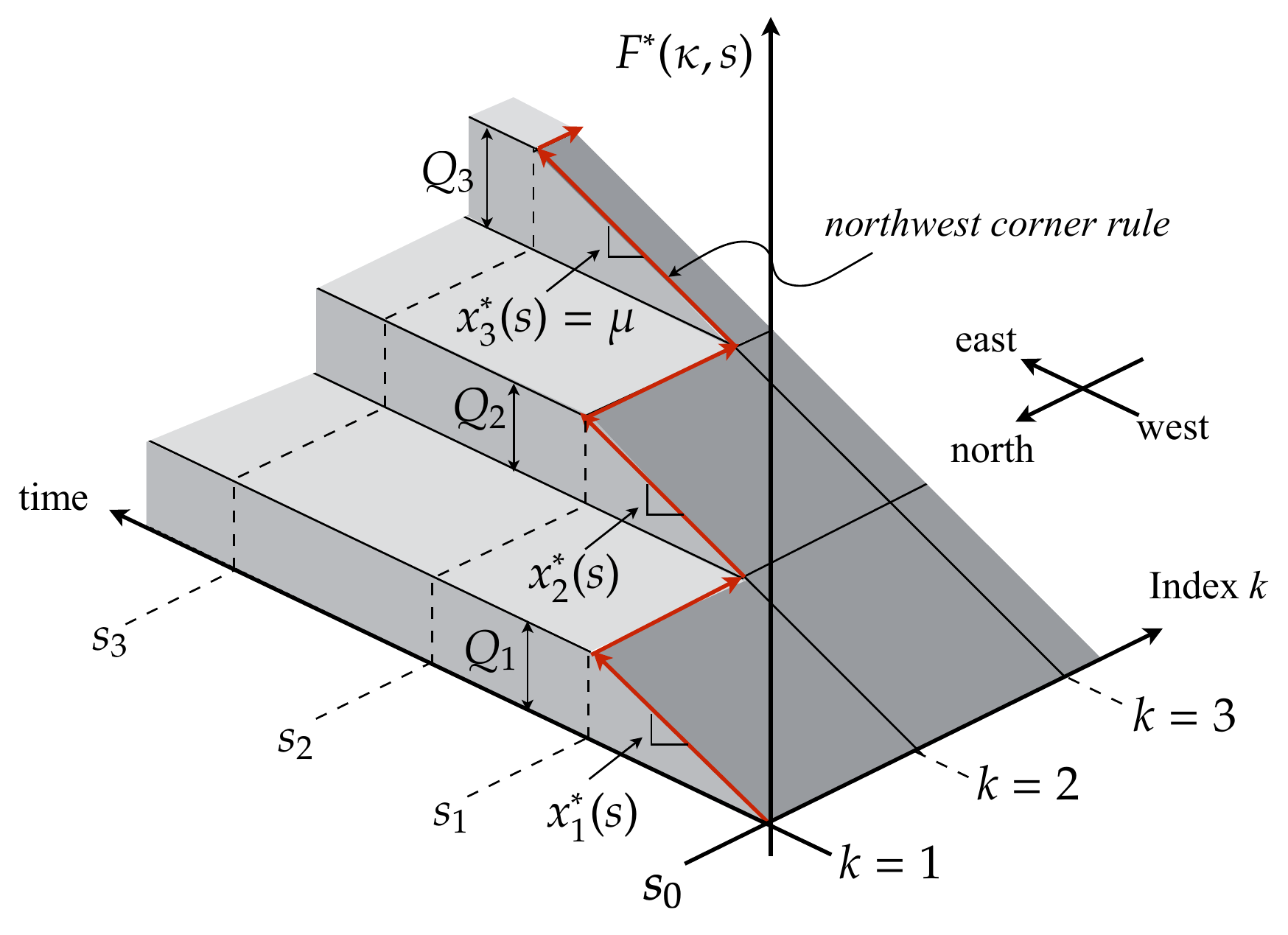}
\caption{Illustration of the solution to the problem \eqref{eq:Cum_OTP-FIFW-obj}}\label{fig:FIFW}
\vspace{-4mm}
\end{center}
\end{figure}

According to Theorem \ref{theo:Frechet-Hoeffding}, we have the following proposition:
\begin{prop}\label{prop:analytical_FIFW} Suppose that Assumption \ref{ass:schedule_FIFW} holds, and user groups are arranged (indexed) in the increasing order of their desired arrival times: $\sigma_{1} < \cdots < \sigma_{K}$.
Then, the solution of [2D-LP(${\bf x}$)] for the equilibrium rush period $\hat{\mathcal{S}}$ is unique and its cumulative form is given by
\begin{align}
	& F^{*}(\kappa, s) = \min \left\{\sum_{k' \leq \kappa}Q_{k'}, \mu \cdot (s - s_{0})\right\} & \forall \kappa \in [1, K], \ \forall s\in \hat{\mathcal{S}} \label{eq:FIFW-analytical}
\end{align}	
\end{prop}

Using this analytical solution \eqref{eq:FIFW-analytical}, we can show the regularity of the equilibrium flow pattern.
Let 
\begin{align}
	s_{k}\equiv s_{0} + \sum_{k' \leq k}Q_{k'}/\mu
\end{align}
be the time when the cumulative supply, $\mu \cdot (s-s_{0})$, is equal to the cumulative demand $\sum_{k' \leq k}Q_{k'}$.
Then, the following equation holds (see also Figure \ref{fig:FIFW}):
\begin{align}
	& F^{*}({k},s_{k}) - F^{*}({k-1}, s_{k-1}) = \int_{s_{k-1}}^{s_{k}}x_{k}^{*}(s) \mathrm{d}s = \mu \cdot (s_{k} - s_{k-1}) = Q_{k} & \forall k = 1, 2, \dots, K.
\end{align}
where $F(0, s_{0})\equiv 0$. 
This equation leads to the following proposition, that is, the FIFW principle.

\begin{prop}
	Suppose that the assumptions in Proposition \ref{prop:analytical_FIFW} hold. Then, the equilibrium flow pattern ${\bf x}^{*}$ has the ``sorting" property, such that all users in group $k$ depart from the bottleneck in a time interval $[s_{k-1}, s_{k}]$ of length $Q_{k}/\mu$ and 
	\begin{align}
		s_{0} < s_{1} < s_{2} < \ \cdots \ < s_{K-1} < s_{K} = s_{0} + T.
	\end{align}
\end{prop}

\begin{figure}[t]
\begin{center}
\includegraphics[width = 70mm,clip]{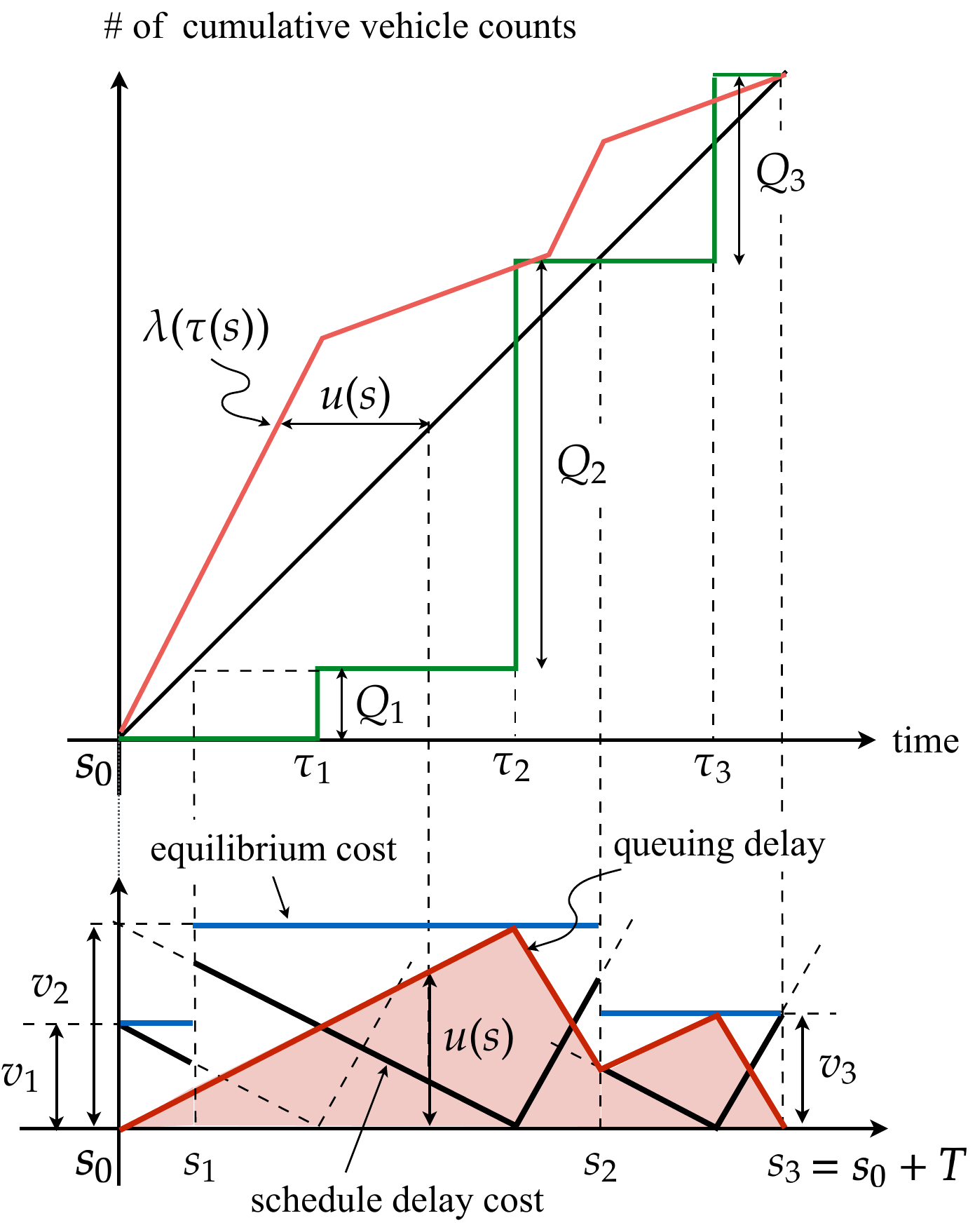}
\caption{Optimal departure flow pattern, cost pattern, and resulting arrival flow pattern ($K = 3$)}\label{fig:FIFW_cum}
\vspace{-4mm}
\end{center}
\end{figure}

The above ``sorting" property also implies that the equilibrium cost pattern ${\bf v}^{*}$ and the associated queuing delays ${\bf u}^{*}$ can be uniquely determined. 
Specifically, the following users' optimal choice condition within each group must hold: 
\begin{align}
	& v_{k}^{*} = u^{*}(s) + c_{k}(s) & \forall s\in [s_{k-1}, s_{k}], \ \forall k\in \mathcal{K}.  \label{eq:cost_FIFW2}
\end{align}
Then we have 
\begin{align}
	& v_{k}^{*} := v_{k+1}^{*} - c_{k+1}(s_{k}) + c_{k}(s_{k}) & \forall k = 1, \dots, K-1. \label{eq:cost_FIFW}
\end{align}
The recursive equation \eqref{eq:cost_FIFW} together with a boundary condition (e.g., the queuing delay for the last user in group $K$ is zero, $v_{K}^{*} = c_{K}(s_{0} + T)$), can be solved easily. 
Note that the Eq.\eqref{eq:cost_FIFW} can also be represented as analytical formulas of ${\bf v}^{*}$ and $s_{0}$ because it consists of simple addition and substitution operations.  
Ultimately, the queuing delay ${\bf u}$ can be obtained analytically by using Eq.\eqref{eq:cost_FIFW2}. 
In addition, as mentioned in the previous section, the cumulative arrival curve can be constructed using the optimal solution $({\bf x}^{*}, {\bf u}^{*})$ above.
Figure \ref{fig:FIFW_cum} shows the relationships between the optimal departure flow pattern, the cost pattern, and the resulting arrival flow pattern for the simplest case (i.e., linear schedule delay cost function).

%%%%%%%%%%%%%%%%%%%%%%%%%%%%%%%%%%%%%%%%%
% 6.  Heterogeneous VOT  %%%%%%%%%%%%%%%%
%%%%%%%%%%%%%%%%%%%%%%%%%%%%%%%%%%%%%%%%%
\section{Analysis of model with heterogeneous value of time}\label{sec:VOT}
In this section, we examine DTCE models to which the optimal transport theory cannot be directly applied.
The problem [2D-LP(${\bf x}$)] is still the equivalent LP problem of the (part of) equilibrium conditions; however, the schedule delay cost function of the models does not satisfy the Monge property. 
Nevertheless, we show that we may be able to solve the equilibrium problem analytically by using an approach that combines the optimal transport theory with a hierarchical decomposition. 
The role of the decomposition in the proposed approach is to generate the lower-level (or sub-) problems that can be analytically solved by invoking the optimal transport theory.
Subsequently, owing to the analytical solutions of the lower-level problems, we may solve the original problem (i.e., the upper-level or master problem) analytically.

As a typical model of the above situation, we consider the DTCE models in which users have the same preferred departure time, $\sigma$, but are classified into $K$ groups (types) differentiated by schedule delay cost functions. 
The schedule delay cost for type $k$ users with departure time, $s$, from the bottleneck is denoted by $c_{k}(\epsilon)$, where $\epsilon \equiv s - \sigma$ is the schedule delay.
The function $c_{k}(\epsilon)$ in this section is assumed to have the following form: 
\begin{subequations} \label{eq:schedule-delay}
\begin{align}
	& c_{k}(s) = \begin{cases}
		\beta_{k} f^{e}(\epsilon) \mspace{15mu} \text{if} \mspace{15mu} \epsilon \leq 0 \mspace{10mu} \text{(early arrival)}\\
		\gamma_{k} f^{l}(\epsilon) \mspace{15mu} \text{if} \mspace{15mu} \epsilon \geq 0 \mspace{10mu} \text{(late arrival)}
	\end{cases}, \\
	& f^{e}(0) = f^{l}(0) = 0
\end{align}	
\end{subequations}
which includes the schedule delay cost functions assumed in conventional models with heterogeneous users \citep[e.g.,][]{Cohen1987,Arnott1988,Arnott1992,Arnott1994a,VandenBerg2011,Ramadurai2010,Liu2015a,Takayama2017} as special cases.
For notational simplicity, the preferred departure time, $\sigma$, will be set to zero hereafter (i.e., $\epsilon = s$).

\begin{figure}[t]
\begin{center}
\includegraphics[width = 70mm,clip]{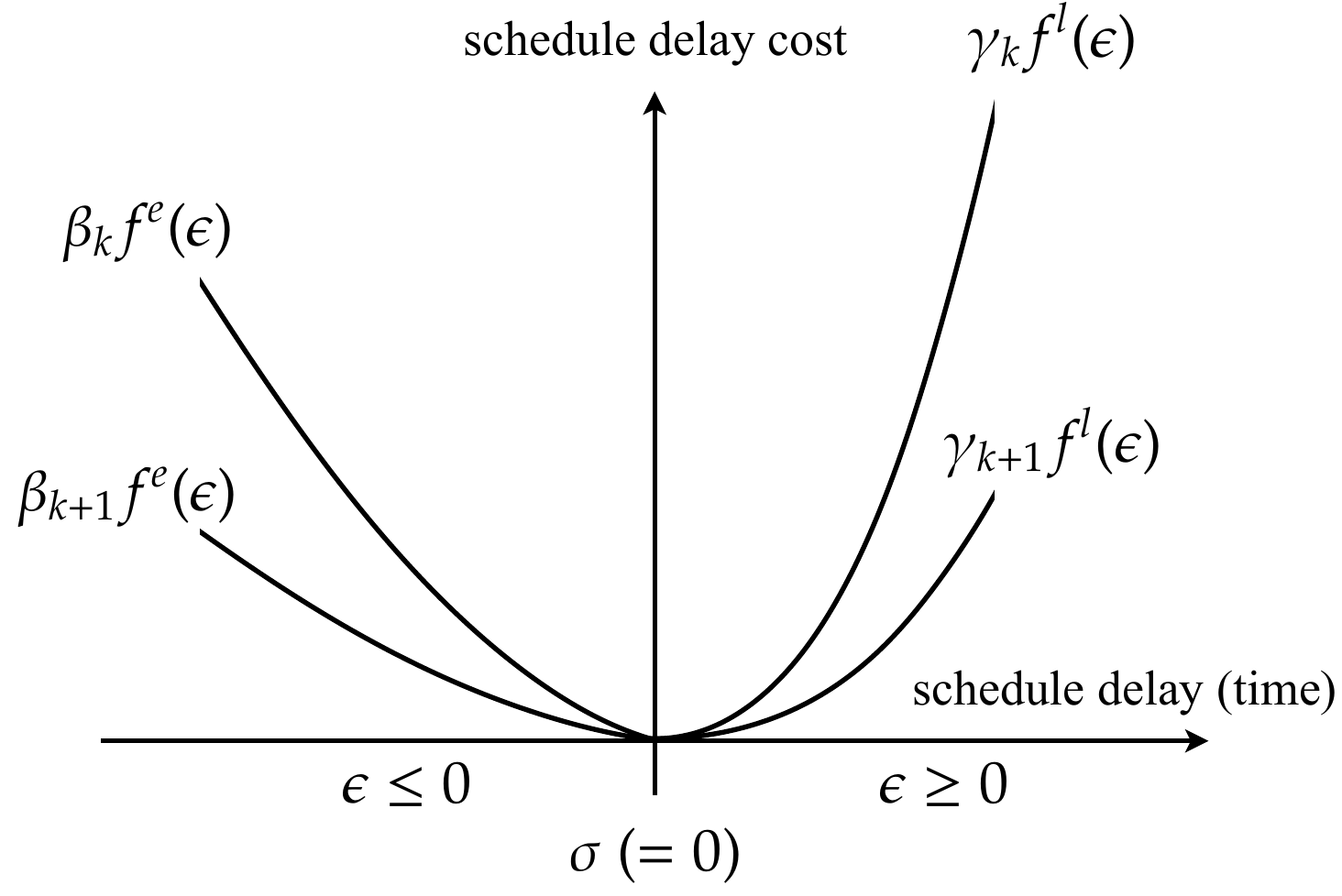}
\caption{Schedule delay $\epsilon$, heterogeneities $k$, and schedule delay cost}\label{fig:ScheduleDelay}
\vspace{-4mm}
\end{center}
\end{figure}

In Sections \ref{sec:Early} and \ref{sec:Late}, we first analyze the models with no late and early arrivals, respectively.
These analyses correspond to solving the lower-level problems of the original problem [2D-LP(${\bf x}$)], which will be formally discussed in the subsequent section. 
In Section \ref{sec:early/late}, we analyze the overall (or uppder-level) problem by exploiting the results in Sections \ref{sec:Early} and \ref{sec:Late}. 
In Section \ref{sec:comparison}, we compare the proposed approach with a semi-analytical approach in the literature.

\subsection{Properties of models with no late arrivals}\label{sec:Early}
We first restrict ourselves to the analysis of models with the following assumption:

\begin{assu}\label{ass:time-indexE}
$f^{e}: \mathbb{R} \rightarrow \mathbb{R}$ is a continuous decreasing function of $\epsilon$ for all $\epsilon \leq 0$ and $f^{e}(0) = 0$, $f^{l}(\epsilon) \rightarrow +\infty$ for all $\epsilon > 0$ (i.e., late arrival is prohibited).
\end{assu}

\begin{lemm}\label{lemm:Early_Monge}
	Suppose that Assumption \ref{ass:time-indexE} holds and that user groups are arranged (indexed) in the decreasing order of their value of time parameters for early arrivals: $\beta_{1} > \beta_{2} > \cdots > \beta_{K}> 0$. Then the function $c_{k}(s)$ defined by \eqref{eq:schedule-delay} satisfies the strict ``inverse" Monge property for $s \leq 0$.
\end{lemm}
\begin{prf} \rm See \ref{sec:proof_l4.1}.     
\end{prf}

As we have seen in Section \ref{sec:DTCandLP}, the DTCE problem can be obtained by solving the LP problem [2D-LP(${\bf x}$)], whose objective is to minimize the total schedule delay cost. 
Because all users have the same preferred departure time $\sigma$, the departure times that are closer to $\sigma$ are chosen at equilibrium. 
Thus, we obtain the following apparent property of an equilibrium departure time window:  
\begin{lemm}\label{lemm:Early_TimeWindows}
	Under Assumption \ref{ass:time-indexE}, the optimal solution ${\bf x}^{*}$ of [2D-LP(${\bf x}$)] satisfies
	\begin{align}
		\begin{cases}
			\sum_{k\in \mathcal{K}}x_{k}^{*}(s) = \mu \mspace{15mu} \text{\rm if} \mspace{15mu} s\in \hat{\mathcal{S}}\subset \mathcal{S}\\
			\sum_{k\in \mathcal{K}}x_{k}^{*}(s) = 0 \mspace{17mu} \text{\rm otherwise} 
		\end{cases}
	\end{align}
	where $\hat{\mathcal{S}}\equiv [-T, 0]$ of length $T= Q/\mu$.
\end{lemm}

These lemmas ensure that [2D-LP(${\bf x}$)] reduces to an optimal transport problem with a strict {\it inverse} Monge property.
However, to use Theorem \ref{theo:Frechet-Hoeffding} directly, it is convenient to formulate the optimal transport problem in its cumulative form with a strict Monge property, as in Subsection \ref{sec:FIFW}.
To do this, let us reverse the time direction: $z \equiv -s$. 
Let $c(\kappa, z)$ and $F(\kappa, z) = \sum_{k' \leq \kappa}\int_{0}^{z}x_{k'}(z)\mathrm{d}z$ be the cost and distribution functions for the time window $\mathcal{Z} = [0, T] \ni z$, respectively. 
Without loss of generality, we also assume that $c(\kappa, z)$ is a continuous and strict submodular function that satisfies $c(k, z) = c_{k}(z), \ \forall k = 1, \dots, K, \ \forall z\in \mathcal{Z}$. 
We then have: 
\begin{align}
& \min_{F\in \mathcal{F}}\int_{[1,K]\times \mathcal{Z}}c(\kappa, z)\mathrm{d}F(\kappa, z)\\
& \text{where} \mspace{15mu} \mathcal{F}\equiv \left\{F(\kappa, T) = {\textstyle \sum_{k' \leq \kappa}}Q_{k'}, \ F(K, t)  = \mu z, \ \forall \kappa \in [1, K], \ \forall z\in \mathcal{Z}\right\}.
\end{align}
We further define the time $s_{k}$ when the cumulative supply $\mu z$ is equal to the cumulative demand $\sum_{k' \leq k}Q_{k'}$, that is, 
\begin{align}
	s_{k} \equiv \sum_{k' \leq k}Q_{k'}/\mu.
\end{align}
By applying Theorem \ref{theo:Frechet-Hoeffding} to this problem, the following proposition is obtained: 
\begin{prop} \label{prop:Early}
	Suppose that Assumption \ref{ass:time-indexE} holds.
	Then, the flow pattern of the departure time choice equilibrium with schedule delay cost function \eqref{eq:schedule-delay} has the following properties:
	\begin{enumerate}
		\item[(1)] The optimal solution of [2D-LP(${\bf x}$)] (i.e., equilibrium flow pattern ${\bf x}^{*}$) is unique and given by the Fr\'echet--Hoeffding distribution.
		\begin{align}
			& F^{*}(\kappa, z) = \min \left\{\sum_{k' \leq \kappa}Q_{k'}, \mu z \right\} & \forall \kappa \in [1, K], \ \forall z\in \mathcal{Z} \label{eq:analyticalEarly}
		\end{align}
		\item[(2)] The equilibrium flow pattern ${\bf x}^{*}$ has the ``sorting" property, such that all users in group $k$ depart from the bottleneck in a time interval $[-s_{k}, -s_{k-1}]$ of length $Q_{k}/\mu$ (i.e., $ s_{k} = s_{k-1} + Q_{k}/\mu$), and
		\begin{align} 
			& -T = -s_{K} < -s_{K-1} < \ \cdots \ < -s_{2} < -s_{1} < -s_{0} \equiv 0.
		\end{align}
	\end{enumerate}
\end{prop}

\begin{figure}[t]
\begin{center}
\includegraphics[width = 80mm,clip]{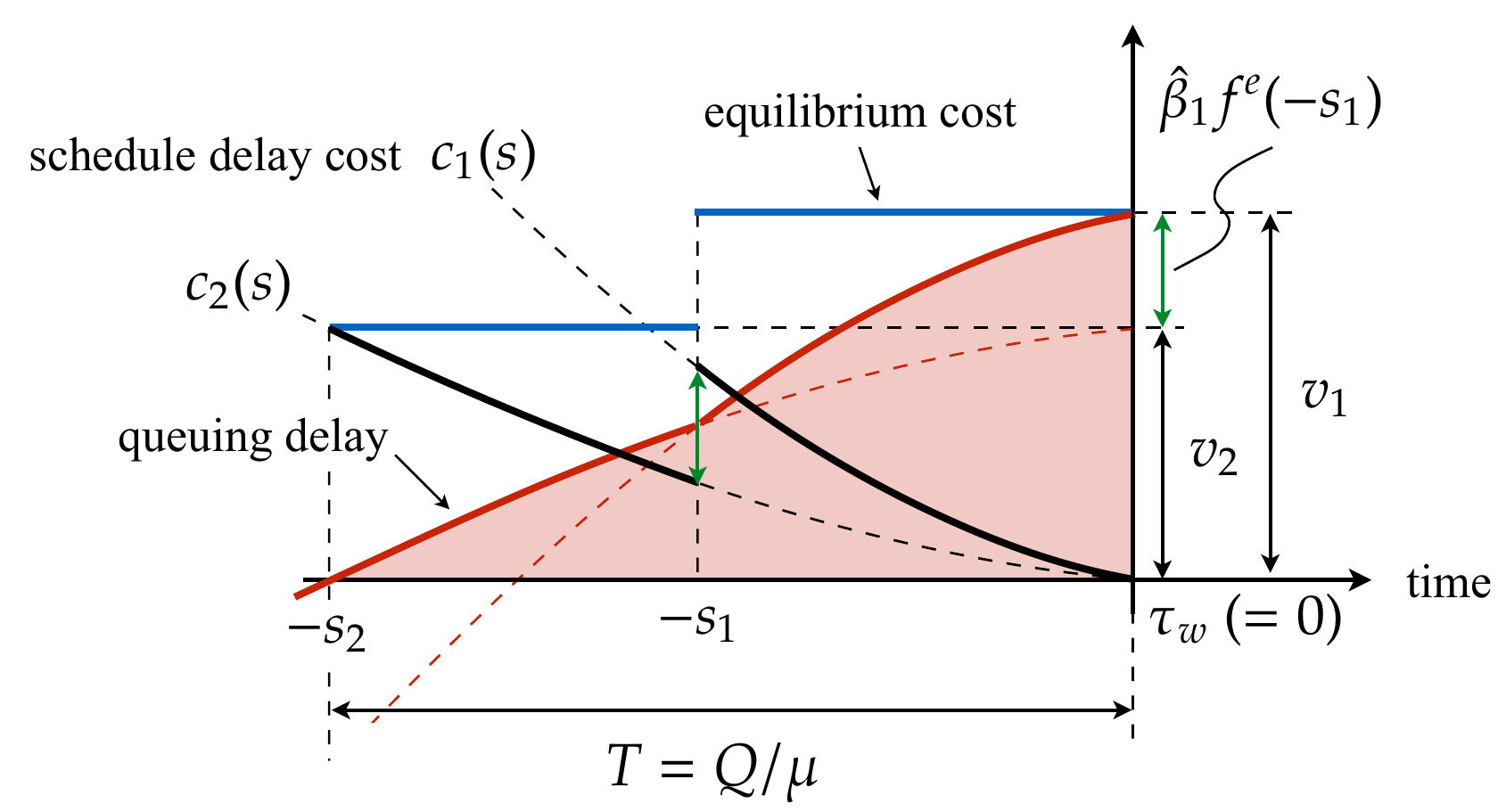}
\caption{Graphical representation of the strong duality for [2D-LP(${\bf x}$)] ($K = 2$, late arrival is prohibited)}\label{fig:CostEarly}
\vspace{-4mm}
\end{center}
\end{figure}

As in the case of Subsection \ref{sec:FIFW}, the equilibrium cost pattern ${\bf v}^{*}$ and associated queuing delays ${\bf u}^{*}$ can be uniquely determined (see Figure \ref{fig:CostEarly}). 
\begin{prop} \label{prop:EarlyCost}
	Suppose that Assumption \ref{ass:time-indexE} holds.
	Then, the cost pattern of the departure time choice equilibrium with schedule delay cost function \eqref{eq:schedule-delay} is unique and is obtained as follows: 
		\begin{align}
			& v_{k}^{*} = \sum_{k' \geq k} \left\{\hat{\beta}_{k'} f^{e}(-s_{k'})\right\} & \forall k \in \mathcal{K} \label{eq:v_early}\\
			& u^{*}(s) = \begin{cases}
				v_{k}^{*} - \beta_{k}f^{e}(s) & \text{\rm if} \ s\in [-s_{k}, -s_{k-1}] \label{eq:u_early}\\
				0 & \text{\rm otherwise}
			\end{cases} & \forall k\in \mathcal{K}
		\end{align}
		where $\hat{\beta}_{k}\equiv \beta_{k} - \beta_{k+1} > 0$ and $\beta_{K+1}\equiv 0$.
\end{prop}
\begin{prf} \rm See \ref{sec:proof_p4.2}.     
\end{prf}

Figure \ref{fig:CostEarly} shows the relationship between the costs at equilibrium.
Specifically, we can see that the following equation holds:  
\begin{align}
	\sum_{k\in \mathcal{K}}(s_{k}-s_{k-1}) v_{k}^{*} - \int_{\mathcal{S}}u^{*}(s)\mathrm{d}s &  = \sum_{k\in \mathcal{K}}\int_{-s_{k}}^{-s_{k-1}}c_{k}(s)\mathrm{d}s \label{eq:Figure5}
\end{align}  
The first term on the LHS of Eq.\eqref{eq:Figure5} shows the total area of the rectangles below the equilibrium trip cost lines (blue bold lines), and the second term shows the red area below the queuing delay curve (red bold line)\footnote{The red dotted lines in Figure \ref{fig:CostEarly} are the so-called {\it isocost queueing curves} in the literature \citep{Lindsey2004}.}. 
The equality holds in equilibrium because the difference in the above-mentioned two areas should be equal to the total area below the schedule delay cost curves (black bold lines).
Furthermore, by recalling that $(s_{k}-s_{k-1}) = Q_{k}/\mu$ and multiplying both sides of Eq.\eqref{eq:Figure5} by $\mu$, we can see the strong duality for [2D-LP(${\bf x}$)].   
\begin{align}
	Z({\bf x}^{*}) =  \sum_{k\in \mathcal{K}}\int_{\mathcal{S}}c_{k}(s)x_{k}^{*}(s)\mathrm{d}s & = \mu\sum_{k\in \mathcal{K}}\int_{-s_{k}}^{-s_{k-1}}c_{k}(s)\mathrm{d}s \notag\\
	& = \sum_{k\in \mathcal{K}}Q_{k} v_{k}^{*} - \mu\int_{\mathcal{S}}u^{*}(s)\mathrm{d}s
	= \check{Z}({\bf u}^{*}, {\bf v}^{*}). \label{eq:strong_duality}
\end{align}

\subsection{Properties of models with no early arrivals}\label{sec:Late}  
In exactly the same manner as we discussed above, we can consider the ``reverse" case, in which the following assumption holds:
\begin{assu}\label{ass:time-indexL}
$f^{l}: \mathbb{R} \rightarrow \mathbb{R}$ is a continuous increasing function of $\epsilon$ for all $\epsilon \geq 0$ and $f^{l}(0) = 0$, $f^{e} \rightarrow +\infty$ for all $\epsilon < 0$ (i.e., early arrival is prohibited).
\end{assu}
\noindent Although this assumption (no early arrivals) may seem strange, the results for this case together with the previous (no late arrivals) case can be used as fundamental building blocks for analyzing more general cases in Subsection \ref{sec:early/late}. 
As in the previous subsection, the following lemma and propositions hold.

\begin{lemm}
	Suppose that Assumption \ref{ass:time-indexL} holds and that user groups are arranged (indexed) in the decreasing order of their value of time parameters for early arrivals: $\gamma_{1} > \gamma_{2} > \cdots > \gamma_{K} > 0$. Then, the function $c_{k}(s)$ defined by \eqref{eq:schedule-delay} satisfies the strict Monge property.
\end{lemm}\vspace{0.5mm}

\begin{prop} \label{prop:Late}
	Suppose that Assumption \ref{ass:time-indexL} holds.
	Then the flow pattern of the departure time choice equilibrium with schedule delay cost function \eqref{eq:schedule-delay} has the following properties:
	\begin{enumerate}
		\item[(1)] The optimal solution of [2D-LP(${\bf x}$)] (i.e., equilibrium flow pattern ${\bf x}^{*}$) is unique and given by the Fr\'echet--Hoeffding distribution.%:
		\begin{align}
			&F^{*}(\kappa, s) = \min \left\{\sum_{k' \leq \kappa}Q_{k'}, \mu s\right\} & \forall \kappa \in [1, K], \ \forall s\in \hat{\mathcal{S}}
		\end{align}
		where $\hat{\mathcal{S}} \equiv [0, T]$ of length $T = Q/\mu$.
		\item[(2)] The equilibrium flow pattern ${\bf x}^{*}$ has the ``sorting" property, such that all users in group $k$ depart from the bottleneck in a time interval $[s_{k-1}, s_{k}]$ of length $Q_{k}/\mu$ (i.e., $s_{k} = s_{k-1} + Q_{k}/\mu$), and
		\begin{align} 
			& s_{0} \equiv 0 < s_{1} < s_{2} < \ \cdots \ < s_{K-1} < s_{K} = T
		\end{align}
	\end{enumerate}
\end{prop}\vspace{0.5mm}

\begin{prop} \label{prop:LateCost}
	Suppose that Assumption \ref{ass:time-indexL} holds.
	Then, the cost pattern of the departure time choice equilibrium with schedule delay cost function \eqref{eq:schedule-delay} is unique and is obtained as follows: 
		\begin{align}
			& v_{k}^{*} = \sum_{k' \geq k} \left\{\hat{\gamma}_{k'} f^{l}(s_{k'})\right\} & \forall k \in \mathcal{K} \label{eq:v_late}\\
			& u^{*}(s) = \begin{cases}
				v_{k}^{*} - \gamma_{k}f^{l}(s) \mspace{15mu} \text{\rm if} \ s\in [s_{k-1}, s_{k}] \label{eq:u_late}\\
				0 \mspace{115mu} \text{\rm otherwise}
			\end{cases} & \forall k\in \mathcal{K}
		\end{align}
		where $\hat{\gamma}_{k}\equiv \gamma_{k} - \gamma_{k+1} > 0$ and $\gamma_{K+1}\equiv 0$.
\end{prop}

\subsection{Properties of models with early and late arrivals}\label{sec:early/late}
At the beginning of Section \ref{sec:VOT}, we briefly mentioned the role of the decomposition in the proposed approach, that is, generating lower-level problems that can be solved analytically.
As also noted, the restricted problems examined in the previous two subsections correspond to such lower-level problems.   
Here, we formally introduce a hierarchical decomposition of the problem [2D-LP(${\bf x}$)], and then show how we can solve the original problem (with a general schedule delay cost function) analytically by exploiting the analytical solutions of the lower-level problems.

Let us consider a general case in which Assumptions \ref{ass:time-indexE} and \ref{ass:time-indexL} are relaxed to allow for both the late and early arrivals.

\begin{assu}\label{ass:time-index}
	$f^{e}: \mathbb{R} \rightarrow \mathbb{R}$ is a continuous decreasing function of $\epsilon$ for all $\epsilon \leq 0$ and $f^{e}(0) = 0$, $f^{l}: \mathbb{R} \rightarrow \mathbb{R}$ is a continuous increasing function of $\epsilon$ for all $\epsilon \geq 0$ and $f^{l}(0) = 0$.
\end{assu}
\noindent This assumption implies that the schedule delay cost function is strictly quasiconvex (or unimodal), which is weaker than Assumption \ref{ass:schedule_FIFW}, under which it is required to be strictly convex.

To decompose the problem [2D-LP(${\bf x}$)], we introduce an auxiliary variable,  $X^{h}_{k}$, representing the total mass of early (if $h = e$) or late (if $h = l$) arrival users in group $k$. 
In addition, let $\mathcal{S}^{h}$ be the time duration $s\leq \sigma$ if $h = e$, $s > \sigma$ if $h = l$. Then, the problem [2D-LP(${\bf x}$)] can be represented as 
\begin{align}
\text{[2D-LP({\bf x}, {\bf X})]} \mspace{15mu} & \min_{{\bf x}\geq {\bf 0}, {\bf X}\geq {\bf 0}}. \ Z({\bf x},{\bf X})\equiv  \sum_{h\in \mathcal{H}}\sum_{k\in \mathcal{K}}\int_{\mathcal{S}^{h}}c_{k}(s) \ x_{k}(s)\mathrm{d}s\\
& \text{subject to} \mspace{15mu} \sum_{k\in \mathcal{K}}x_{k}(s) \leq \mu & \forall s\in \mathcal{S}^{h}, \ h\in \mathcal{H} \label{eq:2DLPxy-con1}\\
& \mspace{80mu}\int_{\mathcal{S}^{h}}x_{k}(s)\mathrm{d}s = X_{k}^{h} & \forall k\in \mathcal{K}, \ h\in \mathcal{H} \label{eq:2DLPxy-con2}\\
& \mspace{82mu}\sum_{h\in \mathcal{H}}X_{k}^{h} = Q_{k} & \forall k\in \mathcal{K} \label{eq:2DLPxy-con3}
\end{align}
where $\mathcal{H} = \{e, l\}$. The equivalence between [2D-LP(${\bf x}$)] and [2D-LP({\bf x}, {\bf X})] is easily confirmed by substituting Eq.\eqref{eq:2DLPxy-con2} into Eq.\eqref{eq:2DLPxy-con3}. 
Following the problem structure in which the variables $\{x_{k}(s)\}$ can be solved if the pattern $\{X_{k}^{h}\}$ is given (but not vice versa), the problem [2D-LP({\bf x}, {\bf X})] can be further viewed as the following equivalent hierarchical problem: 
\begin{align}
	& \min_{{\bf X}\geq {\bf 0}}. \ \sum_{h\in \mathcal{H}}Z_{M}({\bf X}^{h}) \\
	& \text{subject to} \mspace{15mu}  \text{Eq. \eqref{eq:2DLPxy-con3}} \notag\\
	& \text{where} \mspace{15mu} Z_{M}({\bf X}^{h}) \equiv \min_{{\bf x}\geq {\bf 0}}. \ Z_{S}({\bf x} \ | \ {\bf X}^{h})\equiv \sum_{k\in \mathcal{K}}\int_{\mathcal{S}^{h}}c_{k}(s) \ x_{k}(s)\mathrm{d}s \label{eq:LowerP}\\
	& \mspace{180mu}\text{subject to} \mspace{15mu} \text{Eqs. \eqref{eq:2DLPxy-con1} $\ \forall s\in \mathcal{S}^{h}$ and \eqref{eq:2DLPxy-con2}.} \notag
\end{align}

\noindent Because the lower-level (or sub-) problems (for $h= e$ and $h = l$, respectively) are exactly the same problems discussed in the previous subsections, the optimal solution ${\bf x}^{*}$ (for a given ${\bf X}^{h}$) and the associated equilibrium cost (i.e., the optimal Lagrange multiplier) $({\bf u}^{*}, {\bf v}^{h*})$ can be obtained analytically (Propositions \ref{prop:Early}--\ref{prop:LateCost}).
These solutions are useful in the following analysis.

Owing to the strong duality of lower-level problems (i.e., Eq.\eqref{eq:strong_duality}), the upper-level (or master) problem reduces to the following problem: 
\begin{align}
& \min_{{\bf X}\geq {\bf 0}}. \ \sum_{h\in \mathcal{H}}Z_{M}({\bf X}^{h}) = \sum_{h\in \mathcal{H}}\check{Z}({\bf u}^{*}, {\bf v}^{h*}\ | \ {\bf X}^{h})\notag\\ 
& \mspace{140mu} = \sum_{h\in \mathcal{H}}\left[- \mu\int_{\mathcal{S}^{h}}u^{*}(s)\mathrm{d}s + \sum_{k\in \mathcal{K}}X_{k}^{h} v_{k}^{h*}\right] \label{eq:master}\\
& \text{\rm subject to} \mspace{15mu}  \text{\rm Eq. \eqref{eq:2DLPxy-con3}} \notag
\end{align}
In general, an iterative procedure is required to solve the upper-level problem because the dual solution $({\bf u}^{*}, {\bf v}^{h*})$ is not expressed as an explicit function of the variable, ${\bf X}^{h}$.
However, owing to the analytical solution obtained in the previous subsections, we can explicitly transform it into a problem that includes the decision variable ${\bf X}^{h}$.

Recall that there is a one-to-one correspondence between the demand distribution ${\bf X}^{h}$ and the equilibrium departure time vector ${\bf s}^{h}$ (Propositions \ref{prop:Early} and \ref{prop:Late}), that is,  
\begin{align}
	& s_{0}^{h} = 0, \ s_{1}^{h} = X_{1}^{h}/\mu, \ \dots, \ s_{k}^{h} = \sum_{k' \leq k}X_{k'}^{h}/\mu, \ \dots, \ s_{K}^{h} = T^{h} \equiv \sum_{k' \leq K}X_{k'}^{h}/\mu.
\end{align} 
By using these variables and the relationship \eqref{eq:strong_duality}, the upper-level problem can be transformed into the following problem: 
\begin{align}
	& \min_{{\bf s}^{h}}. \ \sum_{h\in \mathcal{H}}Z_{M}({\bf s}^{h})= \mu\sum_{k\in \mathcal{K}}\int_{-s_{k}^{e}}^{-s_{k-1}^{e}}c_{k}(s)\mathrm{d}s + \mu\sum_{k\in \mathcal{K}}\int_{s_{k-1}^{l}}^{s_{k}^{l}}c_{k}(s)\mathrm{d}s \label{eq:CPP-obj}\\
	& \text{subject to} \mspace{15mu}  0\leq s_{1}^{h} \leq \dots \leq s_{k}^{h} \leq \dots \leq s_{K}^{h} \label{eq:CPP-con1} & \forall h\in \mathcal{H}\\
	& \mspace{80mu} s_{k}^{e} + s_{k}^{l} = \sum_{k\leq k'}Q_{k'}/\mu & \forall k\in \mathcal{K} \label{eq:CPP-con2}
\end{align}
It should be noted that this problem is equivalent to [2D-LP({\bf x}, {\bf X})] and [2D-LP({\bf x})], although it can also be interpreted as the upper-level problem of the the hierarchical problem.
As shown in \ref{sec:proof_l4.4}, this problem is a convex programming problem with a strict convex objective function. 
Thus, we obtain the following lemma.

\begin{lemm}\label{lemm:master}Suppose that Assumption \ref{ass:time-index} holds. The optimal solution of problem \eqref{eq:CPP-obj} is unique. 	
\end{lemm}
\begin{prf} \rm See \ref{sec:proof_l4.4}.     
\end{prf}

Finally, we obtain the following proposition.
\begin{prop}\label{prop:early/late}
	Suppose that Assumption \ref{ass:time-index} holds. Then the flow pattern of the departure time choice equilibrium with schedule cost function \eqref{eq:schedule-delay} has the following properties:
	\begin{enumerate}
		\item The optimal solution of [2D-LP(${\bf x}$)] (i.e., equilibrium flow pattern ${\bf x}^{*}$) is unique.
		\item The equilibrium flow pattern ${\bf x}^{*}$ has the ``sorting" property, such that all users in group $k$ depart from the bottleneck in time intervals $[s_{k-1}^{h}, s_{k}^{h}], \forall h\in \mathcal{H}$, and 
		\begin{align}
			-T^{e} = -s_{K}^{e} \leq \ \cdots \ \leq -s_{1}^{e} \leq 0 \leq s_{1}^{l} \leq \ \cdots \ \leq s_{K}^{l} = T^{l}. \label{eq:VOT-sorting}
		\end{align}
		\item If users in all groups depart from the bottleneck both early and late (i.e., strict inequality holds in condition \eqref{eq:VOT-sorting}), the equilibrium values of $s_{k}^{e}$ and $s_{k}^{l}$ are determined by the following equations:
		\begin{align}
			&\begin{cases}
				s_{k}^{l} + s_{k}^{e} = \sum_{k' \leq k}Q_{k'}/\mu  \\
			\hat{\beta}_{k} f^{e}(-s_{k}^{e}) = \hat{\gamma}_{k} f^{l}(s_{k}^{l}) \quad (= v_{k}^{*}-v_{k+1}^{*})
			\end{cases}& \forall k\in \mathcal{K}\label{eq:interior_solution}
		\end{align} 
		where $v_{K+1}^{*} \equiv 0$. 
	\end{enumerate} 
\end{prop}
\begin{prf} \rm See \ref{sec:proof_p4.5}.     
\end{prf}

\begin{figure}[t]
\begin{center}
\includegraphics[width = 72mm,clip]{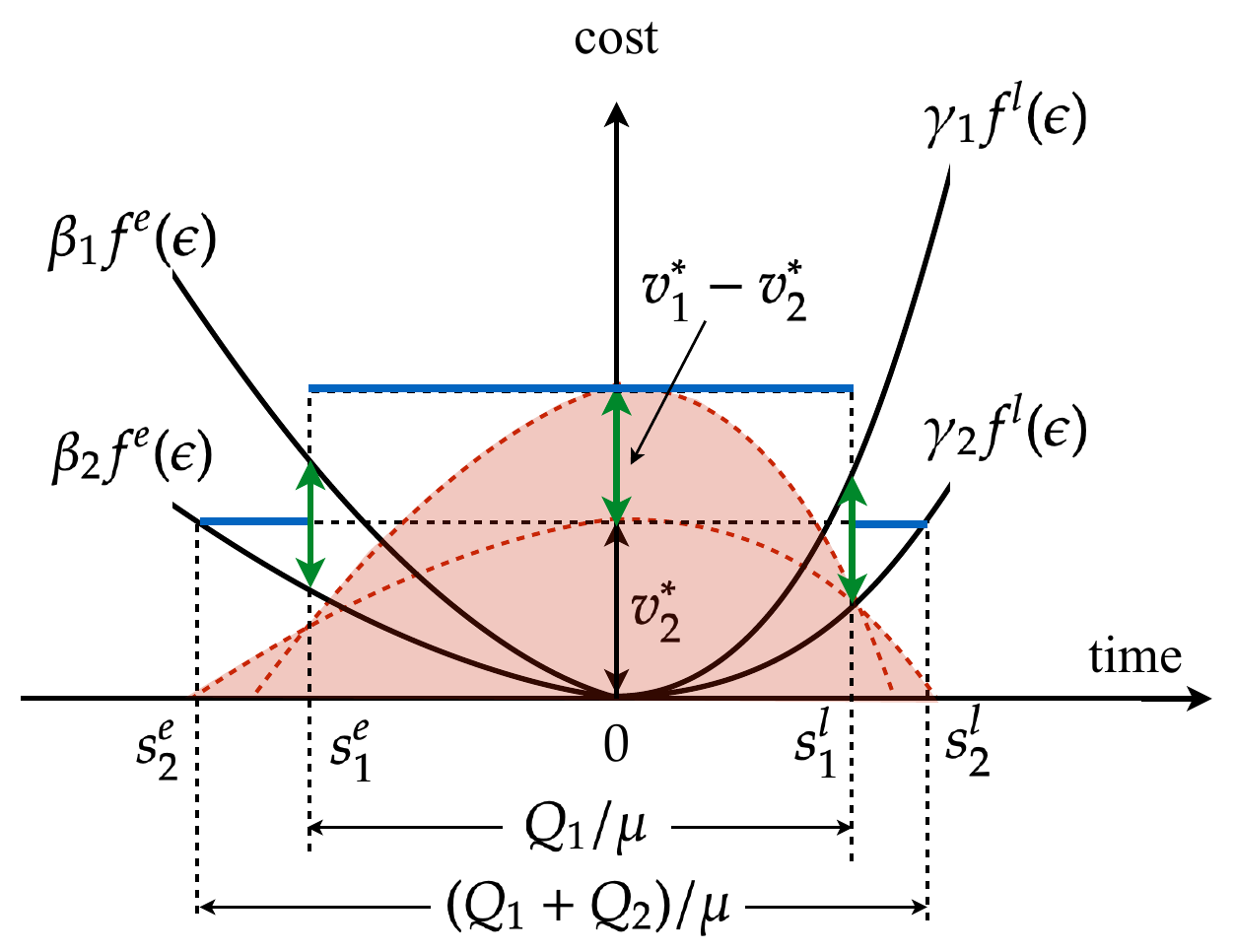}
\caption{Graphical representation of Equation \eqref{eq:interior_solution} ($K = 2$)}\label{fig:EarlyLateb}
\vspace{-4mm}
\end{center}
\end{figure}

\noindent There are a few remarks.
First, the sorting property in this section and the FIFW property seem to be emerging from different mechanisms (i.e., different types of heterogeneities); however, our analyses reveal that both sorting patterns can be explained by the Monge property of the schedule delay cost function in a unified manner. 
Second, property 3 always holds, for example, if $\gamma_{k}/\beta_{k} = \eta$ (constant) for all $k\in \mathcal{K}$, which is a common assumption in the literature  \citep{Vickrey1973,Arnott1988,Arnott1992,Arnott1994a,VandenBerg2011,Takayama2017} (see \ref{sec:proof_p4.5}).
Third, it is evident that Eq.\eqref{eq:interior_solution} for each $k$ is a simple two equations with two variables, that is, $s_{k}^{e}$ and $s_{k}^{l}$ can be easily obtained by independently solving Eq.\eqref{eq:interior_solution} for each $k$ (see Figure \ref{fig:EarlyLateb}). 
If we assume a piecewise linear schedule delay function \cite[e.g.,][]{Arnott1988,Liu2015a}, Eq.\eqref{eq:interior_solution} reduces to elementary linear equations with two variables, which obviously implies that an explicit analytical formula for $s_{k}^{e}$ and $s_{k}^{l}$ is obtained.

\subsection{Relation to an existing semi-analytical approach}\label{sec:comparison}
Before concluding this section, it is worthwhile mentioning a semi-analytical approach proposed by \cite{Liu2015a}.
In their approach, for a given ${\bf X}^{h}$ (i.e., the total mass of early or late arrival users), closed-form trip cost functions are first derived using \cite{Arnott1994a}'s analytical results of the model with a piecewise linear schedule delay function. 
Then, with these trip cost functions (or mappings), they create an equivalent VI problem of the morning commute equilibrium model.  
From the discussions in the previous subsection, the correspondence between \cite{Liu2015a}'s and our proposed approaches are clear: solving the lower-level problems \eqref{eq:LowerP} analytically (Sections \ref{sec:Early} and \ref{sec:Late}) corresponds to a {\it generalized version} of the results of \cite{Arnott1994a}, and the convex programming problem \eqref{eq:CPP-obj}, which is an equivalent problem of the equilibrium model, corresponds to their VI problem. 
The proposed approach is more general than \cite{Liu2015a}'s approach in the sense that we do not impose the restrictive assumption on the schedule delay function.

To discern other differences between the two approaches concretely, let us consider the following VI, which is given as the optimality condition of the upper-level problem \eqref{eq:master}:
\begin{align}
	&\text{Find ${\bf X}^{*}\in  \mathcal{X}\equiv \{{\bf X}\geq {\bf 0} \text{ and Eq.\eqref{eq:2DLPxy-con3}}\}$ such that}  \notag\\
	& \mspace{50mu} \sum_{h\in \mathcal{H}}{\bf v}^{h*}({\bf X}^{h*}) \cdot ({\bf X}^{h} - {\bf X}^{h*}) \geq 0 & {\bf X}\in \mathcal{X}, \label{eq:VIP}
\end{align} 
where the mapping ${\bf v}^{h*}({\bf X}^{h*})$ is obtained by applying the envelope theorem to the dual lower-level problems, i.e., the equilibrium trip cost function:  $\bm{\nabla}Z_{M}({\bf X}^{h}) = \bm{\nabla}Z_{S}({\bf u}^{*}, {\bf v}^{h*} \ | \ {\bf X}^{h}) = {\bf v}^{h*}({\bf X}^{h}), \ \forall h\in \mathcal{H}.$ 
This is almost the same as the VI in \cite{Liu2015a}, except for the measurement unit of trip cost, that is,  the problem \eqref{eq:VIP} is {\it time}-based but that in \cite{Liu2015a} is {\it monetary}-based.
Although the monetary- and time-based formulations can be transformed in to each other by changing the scale of schedule penalty coefficients $\beta_{k}$ and $\gamma_{k}$ for the case of \cite{Liu2015a}, the mathematical properties of the VIs are quite different\footnote{A similar conclusion has been observed in (multiclass) static traffic assignment problems \citep[e.g.,][]{Larsson2002}.}.

More specifically, problem \eqref{eq:VIP} has desirable properties (see \ref{sec:proof_strictMonotone} for the proofs). 
The mapping of problem \eqref{eq:VIP} is strictly monotone; if the schedule delay function is continuously differentiable, the mapping of problem \eqref{eq:VIP} is symmetric and thus a scalar potential exists for the vector field ${\bf v}^{h*}({\bf X}^{h*})$, that is, 
\begin{align}
	\sum_{h\in \mathcal{H}}\oint_{{\bf 0}\rightarrow {\bf X}^{h}}{\bf v}^{h*}(\bm{\omega})\mathrm{d}\bm{\omega} \ \left(= \sum_{h\in \mathcal{H}}Z_{M}({\bf X}^{h})\right).
\end{align}
Therefore, it is easy to obtain the uniqueness of the solutions.
Contrarily, the mapping of \cite{Liu2015a}'s VI is not monotone; it is difficult to establish the uniqueness property of the VI\footnote{
Note that the proof of the uniqueness property in \cite{Liu2015a} is incorrect. 
Specifically, Theorem 2.3 in \cite{More1974}, which was employed in \cite{Liu2015a} for the proof, states that an NCP has at most one solution if the mapping is a $P$-function, i.e., the $P$-function property does not guarantee the uniqueness of VIs when their feasible region is not a {\it positive orthant}; this statement together with the fact that the feasible region of \cite{Liu2015a}'s VI is {\it not} a positive orthant but a simplex means that \cite{Liu2015a}'s proof cannot be valid for their VI. 
}.  
Further, note that the monetary-based VI is not integrable (i.e., it does not reduce to an equivalent LP formulation as does [2D-LP(${\bf x}$)]); thus, the theory of optimal transport and efficient solution methods of LP problems are not applicable to it. 
Note that in \cite{Chen2015a} the time-based version of \cite{Liu2015a}'s VI was studied; it might be reduced to an optimization problem; however they show neither convex programming nor analytical solutions.

In summary, the present approach is a more powerful and promising way to address more general problems.  
As an example, we will analyze an extended problem with more heterogeneities in the next section.

%%%%%%%%%%%%%%%%%%%%%%%%%%%%%%%%%%%%%%%%%
% 7.  Simultaneous Heterogeneity  %%%%%%%
%%%%%%%%%%%%%%%%%%%%%%%%%%%%%%%%%%%%%%%%%
\section{Analysis of models with two types of cost heterogeneity}\label{sec:3D_Monge}
This section demonstrates that the proposed approach is applicable to cases with more general user heterogeneity. 
Specifically, we here consider a generalized DTCE model with two types of cost heterogeneities. 
We then characterize an equilibrium user sorting pattern and obtain an analytical solution by employing both the three-dimensional Monge property and the decomposition of the LP formulation of the model.

\subsection{Review of the $N$-dimensional optimal transport theory}\label{sec:N-OTP}
Before introducing the generalized problem, we briefly show that the Monge property and related useful theorems can be extended to $N$-dimensional transportation problems.
For example, the three-dimensional transportation problem in a discrete space setting is formulated as follows: 
\begin{align}
	\text{[3D-OTP]} \mspace{30mu} & \min_{{\bf x}\geq {\bf 0}}. \ Z_{3D}({\bf x})\equiv  \sum_{i\in \mathcal{I}}\sum_{j\in \mathcal{J}}\sum_{k\in \mathcal{K}}c_{i,j,k}x_{i,j,k} \label{eq:3D-TP_obj}\\
\text{subject to} \mspace{15mu} & \sum_{j\in \mathcal{J}}\sum_{k\in \mathcal{K}}x_{i, j, k} = S_{i} & \forall i\in \mathcal{I} \\
& \sum_{i\in \mathcal{I}}\sum_{k\in \mathcal{K}}x_{i, j, k} = R_{j} & \forall j\in \mathcal{J}\\
& \sum_{j\in \mathcal{J}}\sum_{i\in \mathcal{I}}x_{i, j, k} = Q_{k} & \forall k\in \mathcal{K} 
\end{align}
where ${\bf S}, {\bf R}$, and ${\bf Q}$ are given constants satisfying $\sum_{i\in \mathcal{I}}S_{i} = \sum_{j\in \mathcal{J}}R_{j} = \sum_{k\in \mathcal{K}}Q_{k}$.

The Monge property for higher-dimensional arrays is as follows.
\begin{defi}\rm Let $C$ be an $N$-dimensional array of size $m_{1}\times m_{2}\times \cdots \times m_{N}$. $C$ is termed as a {\it Monge array}, if for all $i_{n} = 1, 2, \dots, m_{n}$, $j_{n} = 1, 2, \dots, m_{n}$, $n = 1, 2, \dots, N$,
\begin{align}
	c[s_{1}, s_{2}, \dots, s_{N}] + c[t_{1}, t_{2}, \dots, t_{N}] \leq c[i_{1}, i_{2}, \dots, i_{N}] + c[j_{1}, j_{2}, \dots, j_{N}]
\end{align}
where $s_{n} \equiv \min (i_{n}, j_{n})$, $t_{n} \equiv \max (i_{n}, j_{n})$ for all $1\leq n\leq N$.
\end{defi}
Note that, as shown in \cite{Aggarwal1988}, an $N$-dimensional array $C$ is a Monge array if and only if, every two-dimensional submatrix is a Monge matrix.
Furthermore, the natural extension of the northwest corner rule solves the problem in a greedy manner without explicit values of cost coefficients if $C$ is a Monge array \citep{Bein1995a}\footnote{Unlike the two-dimensional case, the problem with an {\it inverse} Monge array cannot be solved in a greedy way.}.

As in the two-dimensional case, the continuous Monge property is characterized by a submodular function.
\begin{defi}\rm
A function $f: \mathbb{R}^{N} \rightarrow \mathbb{R}$ is said to satisfy the $N$-dimensional continuous Monge property, if $f$ is submodular with respect to any two of its arguments:
\begin{align}
	& f({\bf x} \land {\bf y}) + f({\bf x} \lor {\bf y}) \leq f({\bf x}) + f({\bf y}) & \forall {\bf x}, {\bf y}\in \mathbb{R}^{N}
\end{align}
where ${\bf x} \land {\bf y}$ and ${\bf x} \lor {\bf y}$ denote the componentwise maximum/minimum of ${\bf x}$ and ${\bf y}$, respectively:
\begin{align}
	&{\bf x} \land {\bf y} \equiv (\max \{x_{1}, y_{1}\}, \max \{x_{2}, y_{2}\}, \dots ,\max \{x_{N}, y_{N}\})  \ \text{\rm for all} \ 1\leq n \leq N  \notag \\
	&{\bf x} \lor {\bf y} \equiv (\min \{x_{1}, y_{1}\}, \min \{x_{2}, y_{2}\}, \dots ,\min \{x_{N}, y_{N}\}) \ \text{\rm for all} \ 1\leq n \leq N. \notag
\end{align}
\end{defi}

Then, \cite{Tchen1980} showed that Theorem \ref{theo:Frechet-Hoeffding} can be generalized for an $N$-dimensional case. 
\begin{theo}\label{theo:Frechet-Hoeffding-N}{\rm \citep{Tchen1980}} Let $F_{1}, \dots, F_{N}$ be distribution functions on $\mathbb{R}$. Furthermore, suppose that $\sup_{F\in \mathcal{F}(F_{1}, \dots, F_{N})} \int_{\mathbb{R}^{N}} c(x_{1},\dots, x_{N})\mathrm{d}F(x_{1},\dots, x_{N}) < \infty$ and a function $c: \mathbb{R}^{N} \rightarrow \mathbb{R}$ satisfies the continuous Monge property.
Then, the Fr\'echet--Hoeffding distribution:
\begin{align}
	& F^{*}(x_{1}, \dots, x_{N}) = \min \{F_{1}(x_{1}), \dots, F_{N}(x_{N})\} & \forall (x_{1}, \dots, x_{N})\in \mathbb{R}^{N} \label{eq:ND-FHdistribution}
\end{align}
is an optimal solution for the following problem:
\begin{align}
	\min_{F\in \mathcal{F}(F_{1}, \dots, F_{N})} \int_{\mathbb{R}^{N}} c(x_{1},\dots, x_{N})\mathrm{d}F(x_{1},\dots, x_{N}).
\end{align}
\end{theo}

\subsection{Simultaneous choice equilibrium model of departure time, location, and job choices}\label{sec:3D-Formulation}
\begin{figure}[t]
\begin{center}
\includegraphics[width = 70mm,clip]{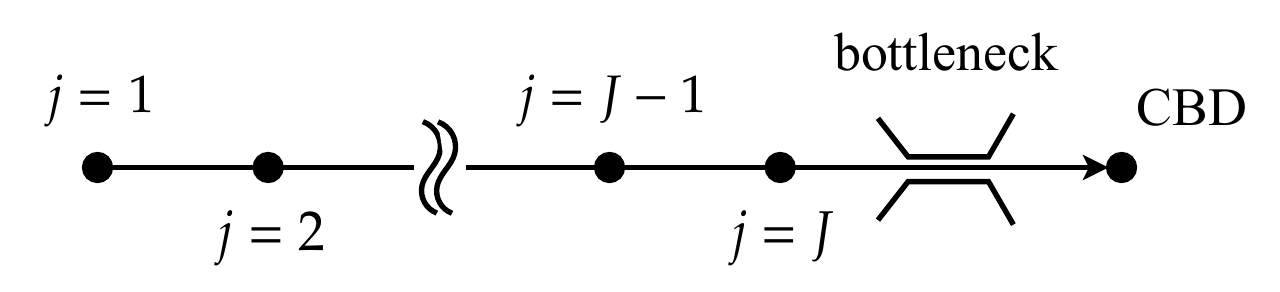}
\caption{Spatial setting inside the city}\label{fig:City}
\vspace{-4mm}
\end{center}
\end{figure}

We show a simultaneous choice equilibrium model of departure time, residential location, and job choices as a generalized DTCE model with two types of cost heterogeneity.
Consider a long narrow city with a spaceless central business district (CBD) where all the job opportunities are located. 
The CBD is located at the edge of the city, and $J$ discrete residential locations are indexed sequentially from the CBD in decreasing order: $J, \dots, 1$ (Figure \ref{fig:City}). 
We denote the set of locations using $\mathcal{J} \equiv \{1, \dots, J\}$. 
The free-flow travel time between location $j$ and the CBD is $l_{j}$, and the maximum allowable number of users who live in each location $j\in \mathcal{J}$ is $R_{j}$. 
There are $K$ types of jobs at the CBD, and the labor demand of each job $k\in \mathcal{K}\equiv \{1, \dots, K\}$ is given by $Q_{k}$.   
The road has a single bottleneck with capacity $\mu$ just upstream of the CBD; thus, all users must pass through this bottleneck to commute.

There are $Q$ ex-ante identical users in the city. 
Each user makes departure time, residential location, and job choices to maximize his/her utility. 
We assume that a user's preference is quasilinear (i.e., the utility function is linear in {\it queuing time units}), and let $w_{k}$, $r_{j}$, and $u(s)$ be the wage of job $k$, land rent at location $j$, and queuing delay at time $s$, respectively. Then, the utility is defined as follows:
\begin{align}
	w_{k} - r_{j} - (c_{j,k}(s) + u(s))	\label{eq:utility}
\end{align}  
where $c_{j,k}(s)$ is a ``commuting disutility" (a generalization of schedule delay cost) with two types of heterogeneity $j$ and $k$, which will be specified later.
Note that all terms in Equation \eqref{eq:utility} are measured in queuing time units.

By combining a user's optimal choice condition and the equilibrium condition in each market, the equilibrium condition is given as follows:
\begin{align}
	&\begin{cases}
		V = w_{k} - r_{j} - u(s) - c_{j,k}(s) \mspace{15mu} \text{if} \mspace{15mu} x_{j,k}(s) > 0\\
		V \geq w_{k} - r_{j} - u(s) - c_{j,k}(s) \mspace{15mu} \text{if} \mspace{15mu} x_{j,k}(s) = 0
	\end{cases} & \forall j \in \mathcal{J}, \ \forall k \in \mathcal{K}, \ \forall s\in \mathcal{S} \label{eq:3D_UserChoice}\\
	& \begin{cases}
		\sum_{j\in \mathcal{J}}\sum_{k\in \mathcal{K}}x_{j,k}(s) = \mu \mspace{15mu} \text{if} \mspace{15mu} u(s) > 0\\
		\sum_{j\in \mathcal{J}}\sum_{k\in \mathcal{K}}x_{j,k}(s) \leq \mu \mspace{15mu} \text{if} \mspace{15mu} u(s) = 0
	\end{cases} & \forall s\in \mathcal{S}\label{eq:q-good}\\
	& \begin{cases}
		\sum_{k\in \mathcal{K}}\int_{\mathcal{S}}x_{j,k}(s)\mathrm{d}s = R_{j} \mspace{15mu} \text{if} \mspace{15mu} r_{j} > 0\\
		\sum_{k\in \mathcal{K}}\int_{\mathcal{S}}x_{j,k}(s)\mathrm{d}s \leq R_{j} \mspace{15mu} \text{if} \mspace{15mu} r_{j} = 0
	\end{cases} & \forall j\in \mathcal{J} \label{eq:J-good}\\
	& \begin{cases}
		\sum_{j\in \mathcal{J}}\int_{\mathcal{S}}x_{j,k}(s)\mathrm{d}s = Q_{k} \mspace{15mu} \text{if} \mspace{15mu} w_{k} > 0\\
		\sum_{j\in \mathcal{J}}\int_{\mathcal{S}}x_{j,k}(s)\mathrm{d}s \geq Q_{k} \mspace{15mu} \text{if} \mspace{15mu} w_{k} = 0
	\end{cases} & \forall k\in \mathcal{K} \label{eq:job-market}\\
	&  \sum_{j\in \mathcal{J}}\sum_{k\in \mathcal{K}}\int_{\mathcal{S}}x_{j,k}(s)\mathrm{d}s = Q \label{eq:general_conservation}
\end{align}
where $V$ represents the equilibrium utility level. 
The first condition \eqref{eq:3D_UserChoice} is the user's optimal choice condition, the second condition \eqref{eq:q-good} is the queuing delay condition, and the third \eqref{eq:J-good} and fourth \eqref{eq:job-market} are the market clearing (or demand-supply equilibrium) conditions for land use and job markets, respectively.

This model is a variant of the DTCE models that include the location choice of users, for example, \cite{Arnott1998a,Gubins2014a,Takayama2017}. 
To keep the model compact, we assume that the wage $w_{k}$ is dependent only on demand-supply conditions \eqref{eq:job-market} and abstract from the effects of the productivity of firms or other realistic determinants of the wage.
The job-location choice modeling in our formulation corresponds to the celebrated Herbert--Stevens model \citep{Herbert1960,Wheaton1974,Berliant2018} in the urban economics literature.
Note that it might be considered that the simultaneous equilibrium model of the different time scale choices is unrealistic, even though their interactions are important factors in relation to appropriate long-term policies \citep[see][for an example]{Osawa2017}. 
However, as shown in the next subsection, we can convert the problem to a bi-level problem where the upper-level (or {\it long-run}) problem is the job-location choice problem and the lower-level (or {\it short-term}) problem is the departure time choice problem. 
As such, the two problems can be analyzed {\it sequentially}.

In parallel with the discussion in the previous sections, this equilibrium model is equivalently written as an LP problem.

\begin{prop} Consider the following infinite-dimensional linear programming problem:  
\begin{align}
	\text{[3D-LP(${\bf x}$)]} \mspace{30mu} & \min_{{\bf x}\geq {\bf 0}}. \ Z({\bf x})\equiv  \sum_{j\in \mathcal{J}}\sum_{k\in \mathcal{K}}\int_{\mathcal{S}}c_{j,k}(s)x_{j,k}(s)\mathrm{d}s \label{eq:3DP-obj}\\
& \text{subject to} \mspace{15mu} \sum_{j\in \mathcal{J}}\sum_{k\in \mathcal{K}}x_{j,k}(s) \leq \mu & \forall s\in \mathcal{S} \label{eq:3DP-con1}\\
&\mspace{80mu} \sum_{k\in \mathcal{K}}\int_{\mathcal{S}}x_{j,k}(s)\mathrm{d}s \leq R_{j} & \forall j\in \mathcal{J} \label{eq:3DP-con2},\\
&\mspace{80mu} \sum_{j\in \mathcal{J}}\int_{\mathcal{S}}x_{j,k}(s)\mathrm{d}s \geq Q_{k} & \forall k\in \mathcal{K} \label{eq:3DP-con3}.
\end{align}
with $\sum_{j\in \mathcal{J}}R_{j} = \sum_{k\in \mathcal{K}}Q_{k} = Q$, and let ${\bf u}$, ${\bf r}$, and ${\bf w}$ be the Lagrange multipliers for the constraints \eqref{eq:3DP-con1}, \eqref{eq:3DP-con2} and \eqref{eq:3DP-con3}, respectively.
Then, the optimal solution $({\bf x}, {\bf u}, {\bf r}, {\bf w})$ of [3D-LP(${\bf x}$)] is consistent with the equilibrium conditions \eqref{eq:3D_UserChoice}--\eqref{eq:general_conservation}. 
\end{prop}

We see that the equivalent optimization problem [3D-LP(${\bf x}$)] of the generalized DTCE problem has almost the same form as the three--dimensional transportation problem [3D-OTP].
Therefore, if the commuting disutility has the Monge property, the equilibrium sorting property would be established by the theory of optimal transport described in Section \ref{sec:N-OTP}.
In the next subsection, we present an example of such disutility functions.

\subsection{Analytical solution and sorting pattern}\label{sec:3D-Solution}
For simplicity of exposition, we assume in the following that all users have a common preferred CBD arrival time, $\sigma = 0$, and ``late arrival" is prohibited \citep[e.g.,][]{Takayama2020a}. 
We then consider the following commuting disutility function: 
\begin{assu}\label{ass:ModelB} The commuting disutility for a type $(j, k)$ user arriving at the destination at time $s$ is given by
\begin{subequations}
\begin{align}
	& c_{j,k}(s) = c_{j, k}^{(1)} + c_{k}^{(2)}(s) + c_{j,k}^{(3)}(s)\\
	& \text{where} \mspace{15mu} c_{j, k}^{(1)} \equiv \alpha_{k}l_{j}, \  c_{k}^{(2)}(s) \equiv \beta_{k}f(s), \ c_{j,k}^{(3)}(s) \equiv \gamma_{k} g(s - l_{j}) 
\end{align}\label{eq:ModelB}
\end{subequations}
where $f(s)$ is a continuous function of ``schedule delay" $s  \ (= s - \sigma)$ satisfying
	\begin{align*}
		\begin{cases}
			f(s) \geq 0 & \text{\rm for} \mspace{15mu} s \leq 0\\
			f(s) \rightarrow \infty & \text{\rm for} \mspace{15mu} s > 0
		\end{cases}, \mspace{15mu}, f'(s) < 0 \mspace{15mu} \forall s \leq 0.
	\end{align*}
	The function $g(t)$ denotes the ``early bird cost"\footnote{This type of disutility can be interpreted as the utility that the user obtains from being at home and has been considered for expressing the heterogeneity in the residential location in the city \citep{Fosgerau2012} and in the trip length within the urban network \citep{Fosgerau2015b}.} for a user departing from his/her home at time $t(s, l_{j}) \equiv s - l_{j} \leq -l_{J}$, which is a continuous function with the following properties:
	\begin{align*}
		g(t(s, l_{j})) > 0, \mspace{15mu} g'(t(s, l_{j})) < 0, \mspace{15mu} g''(t(s, l_{j})) < 0  \mspace{15mu} \forall t \leq -l_{J}
	\end{align*}
	where the coefficients $\alpha_{k}$, $\beta_{k}$, and $\gamma_{k}$ are values of travel time, schedule delay and activity at home, respectively; they are arranged in an increasing order: $\alpha_{1} < \cdots < \alpha_{K}, \ \beta_{1} < \cdots < \beta_{K}, \ \gamma_{1} < \cdots < \gamma_{K}$.
\end{assu}
\noindent The present model can be viewed as a DTCE model with two types of cost heterogeneity: the values of travel time, schedule delay, and activity at home depend on the job $k$\footnote{This is consistent with the assumption in \cite{Takayama2017} of a positive correlation between the VOT of users and their wage.}, and the early bird cost depends on the residential location $j$.
The main difference between the conventional and present models is related to whether the heterogeneity is given exogenously or not. 
Specifically, in a conventional model (with heterogeneous preferred departure times and VOTs) such as \cite{Newell1987} and \cite{Lindsey2004}, a joint distribution of cost heterogeneities is given exogenously. 
In the present model, the marginal distributions of cost heterogeneities, $\{R_{j}\}$ and $\{Q_{k}\}$, are given exogenously; however, the joint distribution $\{X_{j,k}\}$ is determined as a result of job location choice.

\begin{lemm}\label{lemm:3D-Monge2} Suppose that Assumption \ref{ass:ModelB} holds. Then the function $c_{j,k}(s)$ defined by \eqref{eq:ModelB} satisfies the strict 3D-Monge property. 	
\end{lemm}
\begin{prf}	
See \ref{sec:proof_3DMongeB}.
\end{prf}

As in Lemma \ref{lemm:Early_TimeWindows}, combining the fact that [3D-LP(${\bf x}$)] is a total commuting disutility minimization problem with the functional form of $c_{j,k}(s)$ leads to the following property of an equilibrium departure time window. 
\begin{lemm}\label{lemm:3D-Window2} Under Assumption \ref{ass:ModelB}, the optimal solution ${\bf x}^{*}$ of [3D-LP(${\bf x}$)] satisfies
	\begin{align}
		\begin{cases}
			\sum_{j\in \mathcal{J}}\sum_{k\in \mathcal{K}}x_{j,k}^{*}(s) = \mu \mspace{15mu} \text{\rm if} \mspace{15mu} s\in \hat{\mathcal{S}}\subset \mathcal{S}\\
			\sum_{j\in \mathcal{J}}\sum_{k\in \mathcal{K}}x_{j,k}^{*}(s) = 0 \mspace{17mu} \text{\rm otherwise} 
		\end{cases}
	\end{align}
	where $\hat{\mathcal{S}}\equiv [-T, 0]$ of length $T= Q/\mu$. 
\end{lemm}

\noindent Then, the problem [3D-LP(${\bf x}$)] is reduced to a three-dimensional optimal transport problem with a strict 3D-Monge property. That is
\begin{align}
	&\min_{F\in \mathcal{F}} \int_{[1, J]\times [1, K]\times \hat{\mathcal{S}}} c(\zeta, \kappa, s)\mathrm{d}F(\zeta, \kappa, s) \label{eq:3D-COTP}\\
	&\text{where} \mspace{15mu} \mathcal{F}\equiv \big{\{} F(\zeta, K, 0) = {\textstyle \sum_{j'\leq \xi}}R_{j'}, \ F(J, \kappa, 0) = {\textstyle \sum_{k'\leq \kappa}}Q_{k'}, \notag\\
	& \mspace{105mu} F(J, K, s) = \mu \cdot (T + s), \ \forall \xi \in [1, J], \ \forall \kappa \in [1, K], \ \forall s\in \hat{\mathcal{S}}\big{\}}. 
\end{align}
where $c(\zeta, \kappa, s)$ is a continuous and submodular function that satisfies $c(j,k,s) = c_{j,k}(s), \ \forall j = 1, \dots, J, \ \forall k = 1, \dots, K, \ \forall s\in \hat{\mathcal{S}}$.

Based on the theory of optimal transport in Section \ref{sec:N-OTP} (i.e., Theorem \ref{theo:Frechet-Hoeffding-N}), we have the following proposition.

\begin{prop}\label{prop:3D-Monge} Suppose that the commuting disutility, $c_{j,k}(s)$, satisfies Assumption \ref{ass:ModelB}. Then, the solution of the DTCE problem \eqref{eq:3D-COTP} (or [3D-LP(${\bf x}$)]) is unique and given by the Fr\'echet--Hoeffding distribution.
\begin{align}
	F^{*}(\zeta, \kappa, s) = \min\left\{\sum_{j'\leq \zeta}R_{j'}, \sum_{k'\leq \kappa}Q_{k'}, \mu \cdot (T+s)\right\}. \label{eq:analytical-3D}
\end{align}	
\end{prop}

\begin{figure}[t]
\begin{center}
\includegraphics[width = 85mm,clip]{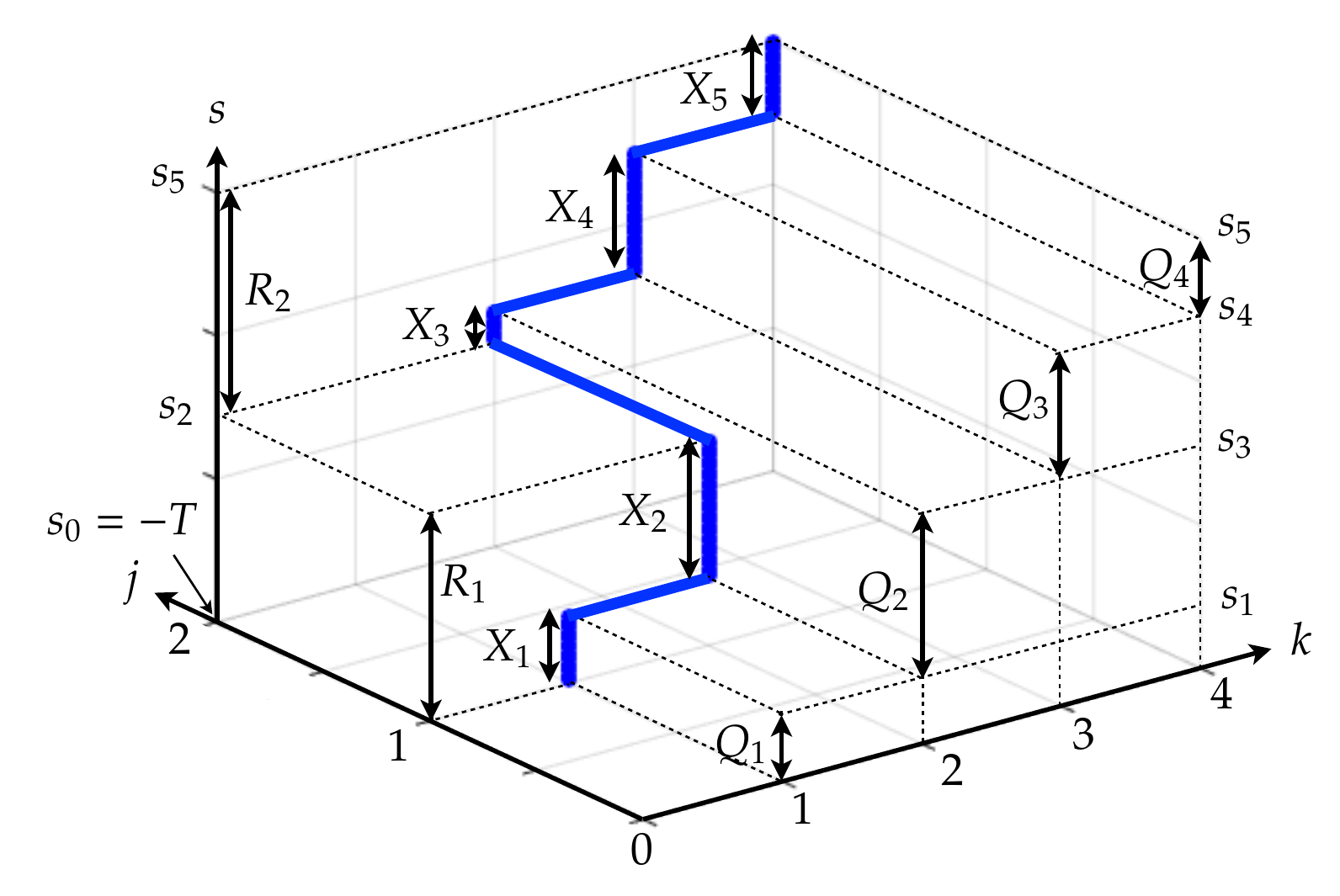}
\caption{Optimal path in the three--dimensional index space}\label{fig:3D}
\vspace{-4mm}
\end{center}
\end{figure}

To interpret the analytical solution, let us introduce a new index $i \in \mathcal{I}$ and a joint distribution $\{X_{i}\}$ of $\{R_{j}\}$ and $\{Q_{k}\}$, which are determined by the northwest corner rule:
\begin{enumerate}
	\setcounter{enumi}{-1}
	\item Initialize the indices with $i:=1$, $j:= 1$, and $k:= 1$.
	\item Set $X_{i} = \min \left\{R_{j}, Q_{k}\right\}$.
	\item Reduce both $R_{j}$ and $Q_{k}$ by $X_{i}$, i.e., $R_{j}:= R_{j} - X_{i}$ and $Q_{k}:= Q_{k} - X_{i}$. If any of $R_{j}$ and $Q_{k}$ become zero, then these indices are increased by one. 
	\item If $\sum_{i'\leq i}X_{i'} = Q$ then $I = i$ and stop; otherwise $i:= i+1$ and go back to Step 2.
\end{enumerate}
This procedure corresponds to solving the minimization problem of the first two terms of the analytical solution \eqref{eq:analytical-3D}. 
By using the unique joint distribution, the analytical solution is reduced to the following form.    
\begin{align}
	F^{*}(\zeta, \kappa, s) =  \min\left\{\sum_{i'\leq i}X_{i'}, \mu \cdot (T+s)\right\} \label{eq:analytical-3D2}
\end{align}
The joint distribution $\{X_{i}\}$ represents the demand distribution of the departure time choice problem. 
As in the previous sections, let us define the time $s_{i}$ so that the cumulative capacity supply $\mu\cdot (T+s)$ is equal to the cumulative demand $\sum_{i'\leq i}X_{i'}$ as follows: 
\begin{align}
	& s_{i} \equiv \sum_{i'\leq i}X_{i'}/\mu - T
\end{align}
Then, the regularities of the equilibrium flow pattern can be summarized as follows.  
\begin{prop}
	The equilibrium flow pattern ${\bf x}^{*}$ has the ``sorting" property such that all users in group $i$ depart from the bottleneck in a time interval $[s_{i-1}, s_{i}]$ of length $X_{i}/\mu$ (i.e., $ s_{i} = s_{i-1} + X_{i}/\mu$), and
		\begin{align}
			& s_{0} \equiv -T < s_{1} < s_{2} < \ \cdots \ < s_{I-1} < s_{I} = 0
		\end{align}
\end{prop}

From the determination process of the joint distribution, we see that job and location groups are sorted in the increasing order of their indices. 
Hence, the proposition means that users who live in locations closer to the CBD and have jobs with higher VOTs depart from the bottleneck at times closer to the desired time. 
Figure \ref{fig:3D} shows an example of the optimal index path of the problem with $J = 2, K = 4$, and $\mu = 1$.
We see that users depart from the bottleneck in the following order along time $s$: $(j, k) = (1, 1), (1, 2), (2, 2), (2, 3), (2, 4)$.

Some final remarks about the analytical solution are presented in order. 
First, the interpretation of the analytical solution provides a clear relationship between short-run and long-run choice problems. 
It is easily seen that the determination of the joint distribution $\{X_{i}\}$ corresponds to the long-run job-location choice problem. 
Given this demand distribution, the problem \eqref{eq:analytical-3D2} determines the short-run departure time choice problem. 
Second, the solution can be understood from the combination of the two 2D-Monge properties: one is with respect to $j$ and $k$, and the other is with respect to $i$ and $s$ (see \ref{sec:3D-decomposition} for details).

%%%%%%%%%%%%%%%%%%%%%%%%%%%%%%%%%%%%%%%%%
% 8.  Concluding remarks and discussion %
%%%%%%%%%%%%%%%%%%%%%%%%%%%%%%%%%%%%%%%%%
\section{Concluding remarks}\label{sec:Conclusion}
This paper presented a systematic approach for analyzing the departure-time choice equilibrium (DTCE) problem of a single bottleneck with heterogeneous commuters. 
The approach is based on the fact that the DTCE is equivalently represented as an LP problem with a special structure, which can be analytically solved by exploiting the theory of optimal transport combined with a hierarchical decomposition. 
Through applying the proposed approach to several types of model with heterogeneous commuters, it is revealed that (i) the essential condition for emerging equilibrium ``sorting patterns," which have been known in the literature, is that the schedule delay functions have the ``Monge property," (ii) equilibrium problems with the Monge property can be solved analytically (or reduced to a very simple system of equations), and (iii) the proposed approach can be applied to a more general problem with more than two types of heterogeneities.

A straightforward extension of this approach would be to analyze models of simultaneous departure time and route choice in a single bottleneck per route network (e.g., a single O-D network with multiple parallel routes) \citep[e.g.,][]{Kuwahara1987,Arnott1992,Iryo2005,Iryo2007,Liu2015a}. As the equilibrium condition can be equivalently represented as a structured LP problem, the approach proposed herein could be applied without significant modification.

Another interesting direction would be to extend the proposed approach to analyze the DTCE problem in a corridor network with multiple bottlenecks \citep[for a specific instance of the corridor problem, see][]{Akamatsu2015a}. 
Although a straightforward formulation of the problem in a corridor network does not reduce to an LP problem, but instead to an LCP, the recent study by \cite{Fu2018} reveals that the solution of the DTCE assignment with homogeneous commuters can be analytically constructed from that of the DSO assignment formulated as a structured LP problem.  Therefore, it is expected that the solution of the DTCE problem with heterogeneous commuters would also be obtained from an LP problem, which can be analyzed by extending the approach presented herein. The partial (positive) answer to this conjecture is provided in \cite{Osawa2017}, and a full treatment is planned for future work.

\appendix

\section{Point queue model}\label{sec:PQ}
According to \cite{Kuwahara1993,Kuwahara1997}, the point queue model can be represented by the following three conditions.
First, the state equation for the number of users queuing at bottleneck $E(t)$ is: 
\begin{subequations}\label{eq:PQ-1}
	\begin{align} 
		& E(t) = A(t) - D(t)\\
		& \mathrm{d}E(t)/\mathrm{d}t = \lambda(t) - \kappa(t)
	\end{align} 
\end{subequations}
where $A(t)$ [$D(t)$] denotes the cumulative arrival [departure] flow at the bottleneck by time, $t$, and their derivatives with respect to time, $t$, are $\lambda(t) = \mathrm{d}A(t)/\mathrm{d}t$ and $\kappa(t) \equiv \mathrm{d}D(t)/\mathrm{d}t$.

The second condition is the exit flow model.
That is, the departure flow rate, $\kappa(t)$, from the bottleneck at time $t$ is 
\begin{align}
	\kappa(t) = \begin{cases}
		\mu & \text{if} \mspace{15mu} E(t) > 0\\
		\min.[\lambda(t), \mu] & \text{if} \mspace{15mu} E(t) = 0 
	\end{cases}. \label{eq:PQ-2}
\end{align}

The final condition concerns the relationship between the state variables and the queuing delay.
In the point queue model, the queuing delay $d(t)$ for a user arriving at the bottleneck at time, $t$, is the horizontal distance between the cumulative arrival curve $A(t)$ and the cumulative departure curve $D(t)$. That is,
\begin{align}
	d(t) = E(t)/\mu \label{eq:PQ-3}
\end{align}
By combining these conditions, \eqref{eq:PQ-1}--\eqref{eq:PQ-3}, we can obtain condition \eqref{eq:OriginalPQ}.

\section{Proof of Lemma 4.1}\label{sec:proof_deltau}
We prove $\Delta \tau(s) = 1 - \Delta b(s) > 0$ (i.e., the non-negative and the finiteness of arrival flow rate) if $\Delta c_{k}(s) > -1$.
Consider first a time window $\mathcal{S}_{k} \subset \mathcal{S}$ when there is at least one group, say $k$, that has a positive departure flow rate, i.e., $x_{k}(s) > 0$. 
Then, the users' optimal choice condition \eqref{eq:UserChoiceB} can be written as 
\begin{align*}
	& c_{k}(s) + b(s) - v_{k} = 0 & \forall s\in \mathcal{S}_{k}.
\end{align*}
This implies that 
\begin{align*}
	& \Delta b(s) = -\Delta c_{k}(s) & \forall s\in \mathcal{S}_{k}.
\end{align*}
Thus, $\Delta c_{k}(s) > -1$ asserts that $\Delta b(s) < 1$. 

For the remaining time window (i.e., $x(s) = 0$), $b(s) = 0$ holds because of the queuing condition \eqref{eq:PQ-capacity}, which implies $\Delta b(s) = 0$. \qed

\section{Strong duality and complementarity slackness for the problem [2D-LP(${\bf x}$)]}\label{sec:Duality}
\subsection{Strong duality}
We first prove the strong duality of 
[2D-LP(${\bf x}$)] and [2D-LP(${\bf u, v}$)],
both of which have infinite-dimensional decision variables.
We first observe the weak duality.
\newtheorem*{theo:B.1}{Theorem B.1}
\begin{theo:B.1}
Let ${\bf x}$ and $({\bf u, v})$
be arbitrary feasible solutions of 
[2D-LP(${\bf x}$)] and [2D-LP(${\bf u, v}$)], respectively.
Then, we have $Z({\bf x})\geq \check{Z}({\bf u}, {\bf v})$.	
\end{theo:B.1}
\begin{prf}
\rm We have
\begin{align}\label{eq:dualitygap}
	Z({\bf x}) - \check{Z}({\bf u}, {\bf v}) & = \sum_{k\in \mathcal{K}}\int_{\mathcal{S}}c_{k}(s)x_{k}(s)\mathrm{d}s + \mu \int_{\mathcal{S}}u(s)\mathrm{d}s - \sum_{k\in \mathcal{K}}Q_{k}v_{k} \notag\\
	& = \sum_{k\in \mathcal{K}}\int_{\mathcal{S}}c_{k}(s)x_{k}(s)\mathrm{d}s + \mu \int_{\mathcal{S}}u(s)\mathrm{d}s - \sum_{k\in \mathcal{K}}v_{k}\int_{\mathcal{S}}x_{k}(s)\mathrm{d}s \notag\\
	& = \int_{\mathcal{S}}\left[\sum_{k\in \mathcal{K}}x_{k}(s)\left(c_{k}(s) + u(s) - v_{k}\right) + u(s)\left(\mu - \sum_{k\in \mathcal{K}}x_{k}(s)\right) \right]\mathrm{d}s\\
	& \geq 0 \notag
\end{align}	
where the inequality is due to \eqref{eq:2DP-con1}, \eqref{eq:2DP-con2}, and \eqref{eq:2DPD-con1}.
We thus have $Z({\bf x})\geq \check{Z}({\bf u}, {\bf v})$. \qed
\end{prf}

Next, we divide the interval $\mathcal{S}$ into several pieces.
Let $N$ be an arbitrary positive integer, 
and let the sub-interval be defined by
\begin{align*}
	\mathcal{S}_{n}&\equiv [\underline{s} + (n-1)\Delta S, \underline{s} + n\Delta S) & n = 1, 2, \dots, N	
\end{align*}
with $\Delta S \equiv (\overline{s} - \underline{s})/N$. Moreover, let $\overline{c}_{k}$ and $\underline{c}_{k}$ be piecewise constant (step/staircase) functions defined as
\begin{align*}
	&\overline{c}_{k} \equiv \begin{cases}
		c_{k,1}^{\max} \mspace{15mu} (s\in \mathcal{S}_{1})\\
		\vdots\\
		c_{k,N}^{\max} \mspace{15mu} (s\in \mathcal{S}_{N})
	\end{cases}, \mspace{45mu}
	\underline{c}_{k} \equiv \begin{cases}
		c_{k,1}^{\min} \mspace{15mu} (s\in \mathcal{S}_{1})\\
		\vdots\\
		c_{k,N}^{\min} \mspace{15mu} (s\in \mathcal{S}_{N})
	\end{cases}	
\intertext{with}
	& c_{k,n}^{\max}\equiv \sup\{c_{k}(s) \ | \ s\in \mathcal{S}_{n}\}, \mspace{45mu} c_{k,n}^{\min}\equiv \inf\{c_{k}(s) \ | \ s\in \mathcal{S}_{n}\}.	
\end{align*}
Then, both functions belong to the following function set 
\begin{align*}
	\mathcal{PC}_{N} \equiv \left\{\text{function} \ f: \mathcal{S}\rightarrow \mathbb{R} \ \big{|} \ f(s) \ \text{is constant over} \ \mathcal{S}_{n} \ \text{for} \ n = 1, \dots, N \right\}	
\end{align*}
Now, consider the following problems:
\begin{itemize}
	\item~[2D-LP(${\bf x}$)-u]: a problem analogous to [2D-LP(${\bf x}$)] with $c_k(s)$ replaced by $\overline{c}_{k}$ and the constraint $x_{k}(s)\in \mathcal{PC}_{N}$ is added. 
	\item~[2D-LP(${\bf x}$)-l]: a problem analogous to [2D-LP(${\bf x}$)] with $c_k(s)$ replaced by $\underline{c}_{k}$ and the constraint $x_{k}(s)\in \mathcal{PC}_{N}$ is added. 
	\item~[2D-LP(${\bf u, v}$)-u]: a problem analogous to [2D-LP(${\bf u, v}$)] with $c_k(s)$ replaced by $\overline{c}_{k}$ and the constraint $u(s)\in \mathcal{PC}_{N}$ is added.
	\item~[2D-LP(${\bf u, v}$)-l]: a problem analogous to [2D-LP(${\bf u, v}$)] with $c_k(s)$ replaced by $\underline{c}_{k}$ and the constraint $u(s)\in \mathcal{PC}_{N}$ is added.
\end{itemize}
Note that these four problems can be regarded as 
LP problems with finite-dimensional decision variables because $x_k(s)$ and $u(s)$ are piecewise constant functions.
Therefore, the strong duality must hold between 
[2D-LP(${\bf x}$)-u] and [2D-LP(${\bf u, v}$)-u], 
and between [2D-LP(${\bf x}$)-l] and [2D-LP(${\bf u, v}$)-l].

Now, let $\rm val[\,\cdot\,]$ and $\mathcal{F}[\,\cdot\,]$ be the optimal value and feasible set of problems $[\,\cdot\,]$. Then, we have the following lemma.

\newtheorem*{lemm:B.1}{Lemma B.1}
\begin{lemm:B.1}\label{lemm:upperlowerlemma}
For any $N>0$, we have \rm
\begin{align*}
	0 \leq {\rm val}\text{[2D-LP(${\bf x}$)]} - {\rm val}\text{[2D-LP(${\bf u, v}$)]} & \leq {\rm val}\text{[2D-LP(${\bf x}$)-u]} - {\rm val}\text{[2D-LP(${\bf u, v}$)-l]}\\
	& = {\rm val}\text{[2D-LP(${\bf x}$)-u]} - {\rm val}\text{[2D-LP(${\bf x}$)-l]}
\end{align*}
\end{lemm:B.1}
\begin{prf}
\rm The first inequality is due to the weak duality, and the equality follows as 
${\rm val}$[2D-LP(${\bf u, v}$)-l] and ${\rm val}$[2D-LP(${\bf x}$)-l] are the primal-dual 
pair of {\it finite-dimensional} LP problems.
Thus, we only show the second inequality.

We first show ${\rm val}\text{[2D-LP(${\bf x}$)]}\leq {\rm val}\text{[2D-LP(${\bf x}$)-u]}$.
Let $\overline Z({\bf x})$ be the objective function of [2D-LP(${\bf x}$)-u].
Then, from $c_k(s)\leq \overline{c}_k(s)$, we have $Z({\bf x})\leq \overline Z({\bf x})$ for any ${\bf x}\geq {\bf 0}$.
Moreover, we have $\mathcal{F}\text{[2D-LP(${\bf x}$)-u]}\subseteq\mathcal{F}\text{[2D-LP(${\bf x}$)]}$ because [2D-LP(${\bf x}$)-u] has additional constraints $x_k(s)\in \mathcal{PC}_{N} (k\in\mathcal{K})$.
Thus, we have ${\rm val}\text{[2D-LP(${\bf x}$)]}\leq {\rm val}\text{[2D-LP(${\bf x}$)-u]}$.

Next, we show ${\rm val}\text{[2D-LP(${\bf u, v}$)]}\geq {\rm val}\text{[2D-LP(${\bf u, v}$)-l]}$.
Notice that $\check Z({\bf u, v})$ serves as the objective function of [2D-LP(${\bf u, v}$)-l] as well as [2D-LP(${\bf u, v}$)]. 
Moreover, we have $\mathcal{F}\text{[2D-LP(${\bf u, v}$)-l]}$ $\subseteq \mathcal{F}\text{[2D-LP(${\bf u, v}$)]}$ since $\mathcal{F}$[2D-LP(${\bf u, v}$)-l] has additional constraints $x_k(s)\in \mathcal{PC}_{N} (k\in\mathcal{K})$ and it follows $\underline{c}_k(s)+u(s)-v_k\geq 0 \ \Rightarrow \ c_k(s)+u(s)-v_k\geq 0$. 
Hence we have ${\rm val}\text{[2D-LP(${\bf u, v}$)]}\geq {\rm val}\text{[2D-LP(${\bf u, v}$)-l]}$. 
This completes the proof. \qed
\end{prf}

Using this lemma, we can show the strong duality 
of [2D-LP(${\bf x}$)] and [2D-LP(${\bf u, v}$)].

\newtheorem*{theo:B.2}{Theorem B.2}
\begin{theo:B.2}
	Suppose that $c_k(s)$ is continuous over $\mathcal{S}$. Then we have \rm ${\rm val}$[2D-LP(${\bf x}$)]$={\rm val}$[2D-LP(${\bf u, v}$)].	
\end{theo:B.2}
\begin{prf}
\rm First, we fix $N>0$ arbitrarily. 
Let $\overline{{\bf x}}^{N}$ and $\underline{{\bf x}}^N$ be the optima of [2D-LP(${\bf x}$)-u] and [2D-LP(${\bf x}$)-l], respectively.	
Noting $\overline{x}^{N}_{k}(s), \underline{x}^N_{k}(s)\in \mathcal{PC}_N$, we denote the $n$-th constant of $\overline{x}^N_k(s)$ and $\underline{x}^N_k(s)$ by $\overline{x}^N_{k,n}$ and $\underline{x}^N_{k,n}$, respectively, for each $k\in \mathcal{K}$.
Moreover, let $\delta^N\geq0$ be defined by
\begin{align*}
	\delta^{N} & \equiv \max\left\{\overline{c}_{k}(s)-\underline{c}_{k}(s) \ \big{|} \ k\in \mathcal{K}, \ s\in \mathcal{S}\right\}	\\
	& = \max\left\{c_{k,n}^{\max}-c_{k,n}^{\min} \ \big{|} \ k\in \mathcal{K}, \ n \in \{1, 2, \dots, N\}\right\}.
\end{align*}
Then, we have 
\begin{align}
	{\rm val}&\text{[2D-LP(${\bf x}$)-u]} - {\rm val}\text{[2D-LP(${\bf x}$)-l]} \notag\\
	& = \sum_{k\in \mathcal{K}}\int_{\mathcal{S}}\overline{c}_{k}(s)\overline{x}_{k}^{N}(s)\mathrm{d}s - \sum_{k\in \mathcal{K}}\int_{\mathcal{S}}\underline{c}_{k}(s)\underline{x}_{k}^{N}(s)\mathrm{d}s \notag\\
	& \leq \sum_{k\in \mathcal{K}}\int_{\mathcal{S}}\overline{c}_{k}(s)\underline{x}_{k}^{N}(s)\mathrm{d}s - \sum_{k\in \mathcal{K}}\int_{\mathcal{S}}\underline{c}_{k}(s)\underline{x}_{k}^{N}(s)\mathrm{d}s \notag\\
	& \leq \sum_{k\in \mathcal{K}}\int_{\mathcal{S}}\delta^{N}\underline{x}_{k}^{N}(s)\mathrm{d}s \notag\\
	& = \delta^{N}\sum_{k\in \mathcal{K}}Q_{k} \label{eq:ineq_deltaN}
\end{align}
where the first inequality follows because $\overline{{\bf x}}^{N}$ is the optimum of [2D-LP(${\bf x}$)-u] and $\underline{{\bf x}}^{N}$ belongs to $\mathcal{F}\text{[2D-LP(${\bf x}$)-u]}$. 
Because $c_k(s)$ is continuous, we have $\lim_{N\to\infty}\delta^N=0$. 
This together with \eqref{eq:ineq_deltaN} and {\bf Lemma B.1} yields ${\rm val}$[2D-LP(${\bf x}$)]$={\rm val}$[2D-LP(${\bf u, v}$)]. This completes the proof. \qed
\end{prf}

In the above theorem, we assumed that $c_k(s)$ is continuous over $\mathcal{S}$. However, even when $c_k(s)$ is discontinuous at a finite number of points on $\mathcal{S}$, we can obtain the same result.
\newtheorem*{coro:B.1}{Corollary B.1}
\begin{coro:B.1}
	Suppose that there exists $\hat s_m\ (m=1,2,\dots,M)$ such that $c_k(s)$ is continuous over $\mathcal{S}\setminus\{\hat s_1,\ldots,\hat s_M\}$. Then we have \rm ${\rm val}$[2D-LP(${\bf x}$)]$={\rm val}$[2D-LP(${\bf u, v}$)].
\end{coro:B.1}
\begin{prf}
\rm	Considering the finitely many sub-intervals 
$[\underline s,\hat s_1),\ [\hat s_1,\hat s_2),\ldots,\ [\hat s_M, \overline s]$, we can prove the corollary in a similar manner.	
\end{prf}

\subsection{Equivalence of complementarity slackness conditions}
Note that \eqref{eq:dualitygap} implies that 
\begin{align*}
	Z({\bf x}) - \check{Z}({\bf u}, {\bf v}) = \int_{\mathcal{S}}\left[\sum_{k\in \mathcal{K}}x_{k}(s)\left(c_{k}(s) + u(s) - v_{k}\right) + u(s)\left(\mu - \sum_{k\in \mathcal{K}}x_{k}(s)\right) \right]\mathrm{d}s.
\end{align*}
Hence, if ${\bf x}$ and $({\bf u, v})$ satisfy \eqref{eq:KT-1}--\eqref{eq:KT-2}, then we have $Z({\bf x})-\check Z({\bf u, v})=0$; that is, ${\bf x}$ and $({\bf u, v})$ are optima of [2D-LP(${\bf x}$)] and [2D-LP(${\bf u, v}$)], respectively.

On the contrary, if ${\bf x}$ and $({\bf u, v})$ optimize [2D-LP(${\bf x}$)] and [2D-LP(${\bf u, v}$)], respectively, then we have $Z({\bf x})-\check Z({\bf u, v})=0$; that is, \eqref{eq:KT-1}--\eqref{eq:KT-2} are satisfied\footnote{In this case, the equalities in \eqref{eq:KT-1} and \eqref{eq:KT-2} hold almost everywhere over $\mathcal{S}$.}.

\section{Proofs of the Monge property}\label{sec:proofs-Monge}
\subsection{Proof of Lemma 5.1}\label{sec:proof_l3.1}
We show that the function $c: \mathcal{K}\times \mathcal{S}\rightarrow \mathbb{R}$ defined by \eqref{eq:schedule_FIFW} is a strict submodular:
\begin{align}
	& c_{k}(s) + c_{k+1}(s') < c_{k+1}(s) + c_{k}(s') 	& \text{\rm for all} \ s < s', \ 1\leq k < K. \label{eq:FIFW_submodular}
\end{align}
Condition \eqref{eq:FIFW_submodular} can be rewritten using $f(\epsilon)$ as follows.
\begin{align}
	\cfrac{f(\epsilon(\sigma_{k+1}, s) + \delta s) - f(\epsilon(\sigma_{k+1}, s))}{\delta s} <  \cfrac{f(\epsilon(\sigma_{k}, s) + \delta s) - f(\epsilon(\sigma_{k}, s))}{\delta s} \label{eq:slope_submodular}
\end{align}
where $\delta s \equiv s' - s > 0$. 
The strict convexity of the function $f$ implies that the slope monotonically increases with the increase in $\epsilon$. 
Combining this property with $\epsilon(\sigma_{k+1}, s) < \epsilon(\sigma_{k}, s)$ ($\because \ \sigma_{k} < \sigma_{k+1}$), we conclude that Eq.\eqref{eq:slope_submodular} and thus Eq.\eqref{eq:FIFW_submodular} hold true.
\qed

\subsection{Proof of Lemma 6.1}\label{sec:proof_l4.1}
We show that the function $c: \mathcal{K}\times \mathcal{S}\rightarrow \mathbb{R}$ defined by \eqref{eq:schedule-delay} is a strict supermodular for $s \leq 0$:
\begin{align}
	& c_{k}(s) + c_{k+1}(s') > c_{k+1}(s) + c_{k}(s') 	& \text{\rm for all} \ s < s'\leq 0, \ 1\leq k < K. \label{eq:Early_supermodular}
\end{align}
The condition \eqref{eq:Early_supermodular} can be rewritten as follows: 
\begin{align}
	(\beta_{k} - \beta_{k+1})\cdot(f^{e}(s)-f^{e}(s'))> 0   
\end{align}
This holds true because $\beta_{k} - \beta_{k+1} > 0$ and $f^{e}(\cdot)$ is a decreasing function of $s$ (or $\epsilon$). \qed

\subsection{Proof of Lemma 7.1}\label{sec:proof_3DMongeB}
We prove that the function $c$: $\mathcal{J}\times \mathcal{K} \times \mathcal{S} \rightarrow \mathbb{R}$ defined by \eqref{eq:ModelB} satisfies the three-dimensional Monge property. To do so, we show that the function $c$ is submodular with respect to any two of its arguments:
\begin{subequations}
\begin{align}
	& c_{j,k}(s) + c_{j,k+1}(s') < c_{j,k+1}(s) + c_{j,k}(s') 	& \text{\rm for all} \ s \leq s' \leq 0, \ 1\leq k < K, \ \forall j\in \mathcal{J}\\
	& c_{j,k}(s) + c_{j+1,k}(s') < c_{j+1,k}(s) + c_{j,k}(s') 	& \text{\rm for all} \ s \leq s' \leq 0, \ 1\leq j < J, \ \forall k\in \mathcal{K}\\
	& c_{j,k}(s) + c_{j+1,k+1}(s) < c_{j+1,k}(s) + c_{j,k+1}(s) 	& \text{\rm for all} \ 1\leq j < J, \ 1\leq k< K, \ \forall s\leq 0
\end{align}\label{eq:3D-MongeB}
\end{subequations}
By a straightforward manipulation, we have the following equivalent conditions:
\begin{subequations}
\begin{align}
	& 	(\beta_{k}-\beta_{k+1})\cfrac{f(s') - f(s)}{\delta s} + (\gamma_{k}-\gamma_{k+1})\cfrac{g(s'-l_{j}) - g(s-l_{j})}{\delta s} > 0\\
	& -\cfrac{1}{\delta l}\left(\cfrac{g(s' - l_{j+1})-g(s - l_{j+1})}{\delta s} - \cfrac{g(s' - l_{j})-g(s - l_{j})}{\delta s}\right) > 0\\
	& -(\alpha_{k} - \alpha_{k+1}) + (\gamma_{k}-\gamma_{k+1})\cfrac{g(s - l_{j+1})-g(s-l_{j})}{\delta l} > 0
\end{align}
\end{subequations}
where $\delta s \equiv s' - s > 0$, and $\delta l \equiv l_{j} - l_{j+1} > 0$. 
The first line holds true because $\beta_{k}-\beta_{k+1}< 0, \gamma_{k}-\gamma_{k+1} < 0$ and $f'(s), g'(t) < 0$; the second line holds true because $g''(t)< 0$; and the third line holds true because $\alpha_{k}-\alpha_{k+1}< 0$. \qed

\section{Other proofs in Section 6}
\subsection{Proof of Proposition 6.2}\label{sec:proof_p4.2}
By solving the following recursive equation:
\begin{align}
	& v_{k}^{*} := v_{k+1}^{*} - c_{k+1}(-s_{k}) + c_{k}(-s_{k}) = v_{k+1}^{*} + \hat{\beta}_{k} f^{e}(-s_{k}) & \forall k = 1, \dots, K-1.
\end{align}
with the boundary condition $v_{K}^{*} = c_{K}(-T) = \hat{\beta}_{K} f^{e}(-s_{K})$, it is evident that the analytical solution \eqref{eq:v_early} is obtained. 
By substituting ${\bf v}^{*}$ into a user's optimal choice condition, we have 
\begin{align}
	& u^{*}(s) = v_{k}^{*} - \beta_{k}f^{e}(s) & \forall s\in [-s_{k}, -s_{k-1}], \ \forall k\in \mathcal{K} 
\end{align}
Combining the queuing condition \eqref{eq:PQ-capacity} with Lemma \ref{lemm:Early_TimeWindows} yields   
\begin{align}
	& \sum_{k\in \mathcal{K}}x_{k}^{*}(s) = 0 < \mu \ \Rightarrow \ u^{*}(s) = 0  &  \forall s\in \mathcal{S}\setminus \hat{\mathcal{S}}.
\end{align}
Thus, we conclude that Eq.\eqref{eq:u_early} holds. 
\qed

\subsection{Proof of Lemma 6.4}\label{sec:proof_l4.4}
Suppose that the objective function of the upper-level problem is convex.
Then the optimality condition of the upper-level problem \eqref{eq:CPP-obj} is given as the following variational inequality (VI) problem.
\begin{align}
	&\text{Find ${\bf s}^{*}\in  \bm{\Omega}\equiv \{\text{Eqs.\eqref{eq:CPP-con1} and \eqref{eq:CPP-con2}}\}$ such that}  \notag\\
	& \mspace{50mu} \sum_{h\in \mathcal{H}} \bm{\nabla} Z_{M}({\bf s}^{h*}) \cdot ({\bf s}^{h} - {\bf s}^{h*}) \geq 0 & {\bf s}^{*}\in  \bm{\Omega}. \label{eq:VIP-s}\end{align} 
where the element of the gradient of the objective function $\bm{\nabla} Z_{M}({\bf s}^{h*})$ is given by
\begin{align}
	\cfrac{\partial Z_{M}({\bf s}^{h})}{\partial s_{k}^{h}} = 
	\begin{cases}
		\mu\hat{\beta}_{k} f^{e}(-s_{k}^{e}) & \text{if} \mspace{15mu} h = e\\
		\mu\hat{\gamma}_{k} f^{l}(s_{k}^{l}) & \text{if} \mspace{15mu} h = l
	\end{cases} 
\end{align}
We then show the strict monotonicity of $\bm{\nabla} Z_{M}({\bf s}^{h*})$. 
For any two non-negative vectors ${\bf s}^{h}$ and $\tilde{{\bf s}}^{h}$ (${\bf s}^{h}\neq \tilde{{\bf s}}^{h}$), we have
\begin{subequations}
\begin{align}
	&(\bm{\nabla} Z_{M}({\bf s}^{e})-\bm{\nabla} Z_{M}(\tilde{{\bf s}}^{e}))\cdot ({\bf s}^{e} - \tilde{{\bf s}}^{e}) = \mu \sum_{k\in \mathcal{K}}\hat{\beta}_{k}\{f^{e}(-s_{k}^{e})-f^{e}(-\tilde{s}_{k}^{e})\}(s_{k}^{e} - \tilde{s}_{k}^{e}) > 0\\
	&(\bm{\nabla} Z_{M}({\bf s}^{l})-\bm{\nabla} Z_{M}(\tilde{{\bf s}}^{l}))\cdot ({\bf s}^{l} - \tilde{{\bf s}}^{l}) = \mu \sum_{k\in \mathcal{K}}\hat{\gamma}_{k}\{f^{l}(s_{k}^{l})-f^{l}(\tilde{s}_{k}^{l})\}(s_{k}^{l} - \tilde{s}_{k}^{l}) > 0
\end{align}	
\end{subequations}
The last inequalities of both equations follow from the fact that $f^{e}(s)$ is a decreasing function of $s$ and $f^{l}(s)$ is an increasing function of $s$.
Thus, we conclude that the upper-level problem \eqref{eq:CPP-obj} is actually a convex programming problem and the optimal solution of it (or the VI problem \eqref{eq:VIP-s}) is unique. \qed

\subsection{Proof of Proposition 6.5}\label{sec:proof_p4.5}
Because properties 1 and 2 of Proposition \ref{prop:early/late} are trivial, we show the proof of the property 3 only. 
If the strict inequality holds in constraint \eqref{eq:CPP-con1} or condition \eqref{eq:VOT-sorting} (i.e., an optimal solution of the problem \eqref{eq:CPP-obj} is in the interior of the feasible region), the equilibrium costs of early and late arrival users in group $k$ must be identical (i.e., $v_{k}^{e*} = v_{k}^{l*}$):
\begin{align}
	& \sum_{k' \geq k} \left\{\hat{\beta}_{k'} f^{e}(-s_{k'}^{e})\right\} = \sum_{k' \geq k} \left\{\hat{\gamma}_{k'} f^{l}(s_{k'}^{l})\right\} & \forall k \in \mathcal{K} \label{eq:interior_condition}
\intertext{which implies}
	& \hat{\beta}_{k} f^{e}(-s_{k}^{e*}) = \hat{\gamma}_{k} f^{l}(s_{k}^{l*}) 
	& \forall k \in \mathcal{K} \label{eq:interior_condition2}
\end{align}
We thus have Eq.\eqref{eq:interior_solution}. \qed

An example is that an interior optimal solution is always feasible and thus exists.
Assume that $\gamma_{k}/\beta_{k} = \eta$ (constant) for all $k\in \mathcal{K}$ and that Eq.\eqref{eq:interior_condition2} is satisfied for all $k\in \mathcal{K}$, that is, 
\begin{align}
	& \begin{cases}
		\mu(s_{k}^{e} + s_{k}^{l}) = \sum_{k'\leq k}Q_{k'}\\
		f^{e}(-s_{k}^{e}) = \eta f^{l}(s_{k}^{l})
	\end{cases} & \forall k\in \mathcal{K} \label{eq:example}
\end{align}
From these equations, we first see that $s_{k}^{e}$ and $s_{k}^{l}$ are positive for $\sum_{k'\leq k}Q_{k'} > 0$.
This is because both $f^{e}(-s)$ and $f^{l}(s)$ are zero at $s = 0$ and are strictly increasing functions of $s > 0$.
We then consider the consecutive groups $k$ and $k+1$ and assume that $s_{k}^{e}\geq s_{k+1}^{e}$. From the second condition of Eq.\eqref{eq:example}, we obtain: 
\begin{align}
	f^{l}(s_{k}^{l}) \geq f^{l}(s_{k+1}^{l}) \mspace{15mu} \Leftrightarrow \mspace{15mu} s_{k}^{l} \geq s_{k+1}^{l}.  
\end{align}
We thus conclude that 
\begin{align}
	s_{k}^{e} + s_{k}^{l} \geq s_{k+1}^{e} + s_{k+1}^{l}.
\end{align}
However, this condition contradicts the first condition of Eq.\eqref{eq:example} (i.e., flow conservation), i.e., $s_{k}^{e} < s_{k+1}^{e}$ must hold. By the same logic, we can obtain $s_{k}^{l} < s_{k+1}^{l}$. 
Hence we finally obtain 
\begin{align}
	& 0 < s_{1}^{h} < s_{2}^{h} < \dots < s_{K-1}^{h} < s_{K}^{h} & \forall h\in \mathcal{H}
\end{align}
which shows the existence of an interior optimal solution.

\subsection{Proof of the strict monotonicity of ${\bf v}^{h*}({\bf X}^{h})$}\label{sec:proof_strictMonotone}
We show that the strict monotonicity of ${\bf v}^{h*}({\bf X}^{h})$ for the case of $h = e$ (we omit the proof for the case of $h = l$ because almost the same discussion holds). 
From the analytical solution \eqref{eq:v_early} and the definition of $s_{k}$, the equilibrium trip cost can be expressed as an explicit function of ${\bf y}^{e}$:
\begin{align}
	& v_{k}^{*e}({\bf X}^{e}) = \sum_{k' \geq k} \hat{\beta}_{k'} \cdot f^{e}({\textstyle -\sum_{m\leq k'}X_{m}^{e}/\mu}) & \forall k \in \mathcal{K}
\intertext{or its vector-matrix form of}
	& {\bf v}^{*}({\bf X}^{e}) = {\bf L}^{\mathrm{T}}{\bf c}(-{\bf s}^{e}({\bf X}^{e})) 
\end{align}
where ${\bf L}$ (${\bf L}^{\mathrm{T}}$) is the lower (upper) triangle matrix, ${\bf c}(-{\bf s}^{e}({\bf X}^{e})) \equiv [\hat{\beta}_{1} \cdot f^{e}(-s_{1}^{e}), \dots, \hat{\beta}_{K} \cdot f^{e}(-s_{K}^{e})]^{\mathrm{T}}$ and ${\bf s}^{e}({\bf X}^{e}) \equiv [X_{1}^{e}/\mu, \dots, \sum_{k' \leq K}X_{k'}^{e}/\mu]^{\mathrm{T}} = (1/\mu){\bf LX}^{e}$.
We thus have, for any two non-negative vectors ${\bf X}^{e}$ and $\tilde{{\bf X}}^{e}$,
\begin{align}
	({\bf v}^{*}({\bf X}^{e})-{\bf v}^{*}(\tilde{{\bf X}}^{e}))\cdot ({\bf X}^{e} - \tilde{{\bf X}}^{e}) & = \{{\bf L}^{\mathrm{T}}({\bf c}(-{\bf s}^{e}({\bf X}^{e})) - {\bf c}(-{\bf s}^{e}(\tilde{{\bf X}}^{e})))\}^{\mathrm{T}} ({\bf X}^{e} - \tilde{{\bf X}}^{e})\notag\\
	& = ({\bf c}(-{\bf s}^{e}({\bf X}^{e})) - {\bf c}(-{\bf s}^{e}(\tilde{{\bf X}}^{e}))^{\mathrm{T}}{\bf L}({\bf X}^{e} - \tilde{{\bf X}}^{e})\notag\\
	& = \mu \sum_{k\in \mathcal{K}}\hat{\beta}_{k}\{f^{e}(-s_{k}^{e})-f^{e}(-\tilde{s}_{k}^{e})\}(s_{k}^{e} - \tilde{s}_{k}^{e}) > 0\notag
\end{align}
The last inequality follows from the fact that $f^{e}(s)$ is a decreasing function of $s$.
Thus, we conclude that the VI problem \eqref{eq:VIP} is a monotone problem. 
Furthermore, if the function $f^{e}(\cdot)$ is continuously differentiable, we have
\begin{align}
	& \cfrac{\partial v_{k}^{*}({\bf X}^{e})}{\partial X_{k'}^{e}} = \begin{cases}
	(-1/\mu)\sum_{k' \geq k} \hat{\beta}_{k'} \cdot f'^{e}(-s_{k'}^{e}) & \text{if} \mspace{15mu} k = k'\\
	(-1/\mu)\sum_{m \geq \max\{k, k'\}} \hat{\beta}_{m} \cdot f'^{e}(-s_{m}^{e}) & \text{if} \mspace{15mu} k \neq k'
	\end{cases}.
\end{align}
This shows that  $\bm{\nabla}{\bf v}^{*}({\bf X}^{e})$ is symmetric and ${\bf v}^{*}({\bf X}^{e})$ is integrable.

\section{Hierarchical decomposition of [3D-LP(${\bf x}$)] }\label{sec:3D-decomposition}
To see the combination of the two Monge properties definitively, we define a new variable $X_{j,k}$ as follows: 
\begin{align}
	& X_{j,k} \equiv \int_{\hat{\mathcal{S}}}x_{j,k}(s)\mathrm{d}s.
\end{align}
Substituting the commuting disutility function \eqref{eq:ModelB} into the objective function of [3D-LP(${\bf x}$)], we can then see that $Z({\bf x})$ can be decomposed into the following two terms:
\begin{align}
& Z({\bf x}) = Z_{M}({\bf X}) + Z_{S}({\bf x})
\intertext{where}
& Z_{M}({\bf X}) \equiv \sum_{j\in \mathcal{J}}\sum_{k\in \mathcal{K}}c_{j,k}^{(1)}X_{j,k}\\
& Z_{S}({\bf x}) \equiv \sum_{j\in \mathcal{J}}\sum_{k\in \mathcal{K}}\int_{\hat{\mathcal{S}}}\left(c_{k}^{(2)}(s) + c_{j,k}^{(3)}(s)\right)x_{j,k}(s)\mathrm{d}s
\end{align}
This implies that [3D-LP(${\bf x}$)] allows the following decomposition:
\begin{align}
	& \min_{{\bf X}\geq {\bf 0}}. \ Z_{M}({\bf X}) + Z_{S}^{*}({\bf X}) \label{eq:3D-LP-decompose-obj}\\
	& \text{subject to} \mspace{15mu} \sum_{k\in \mathcal{K}}X_{j,k} = R_{j} & \forall j\in \mathcal{J} \label{eq:3D-LP-decompose-con2}\\
	& \mspace{85mu} \sum_{j\in \mathcal{J}}X_{j,k} = Q_{k} & \forall k\in \mathcal{K} \label{eq:3D-LP-decompose-con3}\\
	& \text{where} \mspace{15mu} Z_{S}^{*}({\bf X})\equiv \min_{{\bf x}\geq {\bf 0}}. \ Z_{S}({\bf x} \ | \ {\bf X}) \label{eq:3D-sub}\\
	& \mspace{145mu}\text{subject to} \mspace{15mu} \sum_{j\in \mathcal{J}}\sum_{k\in \mathcal{K}}x_{j,k}(s) = \mu & \forall s\in \hat{\mathcal{S}} \label{eq:3D-sub-con1}\\
	& \mspace{225mu} \int_{\hat{\mathcal{S}}}x_{j,k}(s)\mathrm{d}s = X_{j,k} & \forall j\in \mathcal{J}, \ \forall k\in \mathcal{K} \label{eq:3D-sub-con2}
\end{align}
The lower-level problem (or subproblem) determines the short-run DTCE for a given joint distribution of $\{R_{j}\}$ and $\{Q_{k}\}$. 
The upper-level (or master) problem determines the long-run job-location choice equilibrium based on the cost including the short-run effect $Z_{S}^{*}({\bf X})$. 
While the upper-level and sub-problems should be solved iteratively in general, the analytical solution implies that these problems can be solved sequentially. 
Specifically, the upper-level problem can be solved by the northwest corner rule (i.e., the objective function \eqref{eq:3D-LP-decompose-obj} satisfies the 2D-Monge property with respect to $j$ and $k$), and the lower-level problem can also be solved in a similar manner for a given optimal solution $\{X_{j,k}\}$ or equivalently $\{X_{i}\}$ (i.e., the objective function \eqref{eq:3D-sub} satisfies the 2D-Monge property with respect to $s$ and $i$).

\bibliographystyle{elsarticle-harv}
\bibliography{DTCandOTP}

%% not want to use model2-names.bst.



%% \bibitem must have one of the following forms:
%%   \bibitem[Jones et al.(1990)]{key}...
%%   \bibitem[Jones et al.(1990)Jones, Baker, and Williams]{key}...
%%   \bibitem[Jones et al., 1990]{key}...
%%   \bibitem[\protect\citeauthoryear{Jones, Baker, and Williams}{Jones
%%       et al.}{1990}]{key}...
%%   \bibitem[\protect\citeauthoryear{Jones et al.}{1990}]{key}...
%%   \bibitem[\protect\astroncite{Jones et al.}{1990}]{key}...
%%   \bibitem[\protect\citename{Jones et al., }1990]{key}...
%%   \harvarditem[Jones et al.]{Jones, Baker, and Williams}{1990}{key}...
%%


% \end{thebibliography}

\end{document}